\documentclass[twocolumn]{autart}    % Enable this line and disable the 
\usepackage{graphicx}          
\usepackage{url}
\usepackage{color}
\usepackage{amsfonts}
\usepackage{cuted} 
\usepackage{hyperref}
\usepackage{amssymb}
\usepackage{bbold}
\usepackage{tikz}
\usepackage{siunitx}
\usepackage{tabularray}
\usepackage{tabularx}
\usepackage{multirow}

\usepackage{physics}

\newtheorem{theorem}{Theorem}[section]
\newtheorem{lemma}[theorem]{Lemma}
\newtheorem{proposition}[theorem]{Proposition}
\newtheorem{remark}[theorem]{Remark}

\newtheorem{definition}[theorem]{Definition}

\newtheorem{example}{Example}
\newtheorem{assumption}[theorem]{Assumption}

\usepackage{subcaption}

\usepackage{amsmath}
\usepackage{algorithm}
\usepackage{algorithmic}

\usepackage{mathtools}
\usepackage{fancyhdr}
\usepackage[T1]{fontenc}
\usepackage{siunitx}
\usepackage{multirow}
\usepackage{placeins}
\usepackage{booktabs}
\usepackage{rotating}
\usepackage{array}
\usepackage{dsfont}
\usepackage{bbold}
\usepackage{adjustbox}
\usetikzlibrary{calc}
\usepackage{xcolor}
\definecolor{niceGreen}{rgb}{0.1, 0.625, 0.1}
\usepackage{fdsymbol}
\usepackage{derivative}
\usepackage[breakable]{tcolorbox}

\usepackage[normalem]{ulem}

\newcommand*{\aleq}{\stackrel{(\text{a})}{\leq}}

\newcommand*{\beq}{\stackrel{(\text{b})}{=}}
\newcommand*{\ceq}{\stackrel{(\text{c})}{=}}

\newcommand*{\ageq}{\stackrel{(\text{a})}{\geq}}



\usetikzlibrary{positioning, fit, calc}   
\tikzset{block/.style={draw, thick, text width=0.8\columnwidth ,minimum height=0.5cm, align=center},   
line/.style={-latex}     
}  

\usepackage{mathalfa}

\definecolor{dartmouthgreen}{rgb}{0.05, 0.5, 0.06}
\begin{document}
\begin{frontmatter}

\title{Synthesis of safety certificates for discrete-time uncertain systems via convex optimization\thanksref{footnoteinfo}} % Title, preferably not more 
     
\thanks[footnoteinfo]{Marta Fochesato's research visit at the University of Oxford was supported by the ETH Z\"urich Doc.Mobility Fellowship Programme.}
\author[eth]{Marta Fochesato}\ead{mfochesato@ethz.ch}, 
\author[oxford]{Han Wang}\ead{han.wang@eng.ox.ac.uk}, 
\author[oxford]{Antonis Papachristodoulou}\ead{antonis@eng.ox.ac.uk}, 
\author[oxford]{Paul Goulart}\ead{paul.goulart@eng.ox.ac.uk},

\address[eth]{Automatic Control Laboratory (IfA), ETH Z\"{u}rich, Z\"{u}rich, Switzerland}  
\address[oxford]{Department of
Engineering Science, University of Oxford, OX1 3PJ, Oxford, UK}

                             % keyword list or with the 
                                          % help of the Automatica 
                                          % keyword wizard

\begin{abstract}                          % Abstract of not more than 200 words.
We study the problem of co-designing control barrier functions and linear state feedback controllers for discrete-time linear systems affected by additive disturbances. For disturbances of bounded magnitude, we provide a semi-definite program whose feasibility implies the existence of a control law and a certificate ensuring safety in the infinite horizon with respect to the worst-case disturbance realization in the uncertainty set. For disturbances with unbounded support, we rely on martingale theory to derive a second semi-definite program whose feasibility provides probabilistic safety guarantees holding joint-in-time over a finite time horizon. We examine several extensions, including (i) encoding of different types of input constraints, (ii) robustification against distributional ambiguity around the true distribution, (iii) design of safety filters, and (iv) extension to general safety specifications such as obstacle avoidance.
\end{abstract}
\end{frontmatter}

\section{Introduction}

Safety requirements are a critical concern in control theory, with many real-world applications calling for certificates ensuring that the system avoids entering unexpected regions, or remains inside a prescribed safe set during the control task length despite the presence of uncertainties. Examples are ubiquitous, ranging from air traffic control \cite{prandini2008application} to autonomous vehicles \cite{bojarski2016end} and robotics \cite{reher2021dynamic}. Consequently, in the wake of the ongoing emergence of ever-more complex safety-critical control problems, there is a growing need for optimal control tools that can provide reliable certificates under uncertainty. 

Consider a stochastic dynamical system $x_{t+1} = f(x_t, u_t, w_t)$ with state $x \in \mathcal{X} \subseteq \mathbb{R}^n$, input $u \in \mathcal{U} \subseteq \mathbb{R}^m$, and a stochastic disturbance $w_t \in \mathcal{W} \subseteq \mathbb{R}^d$. Concurrently, consider an (application-specific) safe set 
$\mathcal{S} \subseteq \mathcal{X}$ describing the portion of the state space in which the system state is allowed to lie. For example, $\mathcal{S}$ may describe the road boundaries in autonomous driving applications, or the allowable voltage and current ranges in power converter operations. Two fundamental problems arise in the safety enforcement problem: (i) for a set of initial conditions $\mathcal{I} \subseteq \mathcal{S}$, \textit{verify} whether there exists a feedback control law $u = \pi(x)$ such that the trajectories starting from $\mathcal{I}$ are safe, i.e., they remain in $\mathcal{S}$ along a prescribed time interval, and (ii) \textit{design} such a safe feedback controller in a reliable and computationally efficient manner. 

The discrete-time Control Barrier Function (CBF) approach answers both of these questions via a continuous function that satisfies certain properties. Roughly speaking, a CBF aims to separate the safe and unsafe regions of the state space by its zero level set and ultimately enables controller synthesis for safety requirements specified by forward invariance of a set using a Lyapunov-like condition. Despite practical success in many applications, such as Segway balancing, bipedal walking robots, and quadrotors (see, for example,  \cite{ames2014control} and references therein), the CBF approach suffers from a fundamental challenge as the design of viable and non-conservative CBF candidates is non-trivial even for linear deterministic systems. In fact, even verifying that a candidate function is a CBF is NP-hard \cite{clark2022semi}.

\subsection{Related works}
While the continuous-time setting has undoubtedly received more attention, several methods have been proposed in recent years to synthesize control barrier functions for discrete-time systems. For example, reachability-based methods \cite{bertsekas1972infinite} solve an optimal control problem returning a set of states starting from which the dynamical system is certified to stay within a prescribed safe set during a given time interval. Given this formulation, one can show that the value function of the optimal control problem is a time-varying CBF \cite{massiani2022safe}. However, solving the optimal control problem is generally challenging as it involves the solution of a Hamilton-Jacobi Partial Differential Equation \cite{margellos2011hamilton,mitchell2005time} and requires ad-hoc numerical heuristics. 

Learning-based methods parametrize both the control law and the CBF, and rely on the availability of a finite dataset to determine the optimal parameters minimizing a suitable loss function \cite{srinivasan2020synthesis,robey2020learning,dawson2023safe}. For example, \cite{wang2024simultaneous,long2024distributionally} consider a neural network-based parametrization. This class of methods can address highly nonlinear and high-order systems; however, they lack rigorous guarantees and require large amounts of training data that are not easy to collect in general \cite{yang2023synthesizing}.

Finally, we recall optimization-based approaches, mostly in the form of Sum-Of-Squares (SOS) optimization. For polynomial dynamics and semi-algebraic safe sets, SOS turns the algebraic conditions for
CBF into polynomial positivity conditions, and casts them using SOS hierarchies \cite{prajna2006barrier,wang2023assessing}. Compared to the previous methods, SOS allows for rigorous safety guarantees and efficient computation (unlike learning-based methods), and tractable solutions (unlike reachability-based methods). However, when
controllers need to be simultaneously designed, a commonly used method alternates between
synthesizing the controller and the CBF in sequential SOS
programs \cite{schneeberger2024advanced,wang2023safety}. Unfortunately, alternating methods suffer from two major drawbacks: first, a feasible initial condition must be provided (which may in general be nontrivial to find \cite{zhao2023convex}, \cite{schneeberger2024advanced}); second, no convergence guarantees are provided \cite{wang2023assessing}. Finally, exploiting SOS programming duality, moment problems based
on occupation measures can be formulated to solve the address the synthesis problem \cite{korda2014convex}. 

Whenever the system dynamics are affected by uncertainty, the design procedure faces additional challenges. Recent literature has attempted to address this additional complexity. The seminal work in \cite{prajna2007framework} only considers the worst-case and stochastic safety verification problems, hence disregarding the control design task. More recently, \cite{jagtap2020formal} addresses the control design task for stochastic discrete-time systems subject to temporal logic specifications. Further, \cite{singletary2022safe} introduces so-called risk control barrier functions,  which are compositions of barrier functions and dynamic coherent risk measures, to enable a risk-aware safety analysis. However, both papers assume a candidate CBF to be known \emph{a priori} to the control designer. As a result, the convex co-design of CBFs and feedback controllers for discrete-time stochastic systems has not found a comprehensive treatment in the literature and is still an open problem. 
% Additionally, \cite{clark2021control} outlines a framework
% for CBF design in continuous-time stochastic systems in the presence of Gaussian process and measurement noise. However, the soundness of the results was recently questioned in \cite{so2025comment} and partly addressed in \cite{so2023almost} under a more restrictive setting.

\subsection{Contributions and Outline}
 The main goal of this manuscript is to address this challenge by providing a convex design procedure for the joint synthesis of a CBF and a feedback controller for the class of discrete-time uncertain linear systems affected by additive stochastic noise. We begin by describing in Section \ref{sec:2} meaningful notions of safety for uncertain systems based on the properties of the disturbance vector. In particular, we distinguish between infinite-horizon safety when the support $\mathcal{W}$ is bounded; and finite-horizon safety that holds joint-in-time over a finite time horizon when $\mathcal{W}$ is unbounded. Next, in Section \ref{sec:main} we formalize our main contribution and show that the above safety requirements can be encoded via convex programs, hence bypassing the need for alternating methods. Further, in Section \ref{sec:4} we consider several extensions: most notably, we consider the inclusion of different types of input constraints, the robustification against ambiguity in the probability distribution governing the disturbance, the design of safety filters, and the extension to general safety specifications. Finally, we show the effectiveness of our novel design procedure in simulation (Section \ref{sec:5}). Section \ref{sec:6} concludes the paper.

\subsection{Notation and Basic Preliminaries}
$\mathbb{R}^n, \mathbb{R}^n_{\geq 0}, \mathbb{N}$ denote the space of $n-$dimensional real numbers, $n-$dimensional nonnegative numbers and natural numbers, respectively. 

Let $\mathcal{P}(\Omega)$ denote the power set and $\mathbf{B}(\Omega)$ denote the Borel field of some set $\Omega$. A set ${F} \subset \mathbb{R}^n$ is Borel measurable if ${F} \in \mathbf{B}(\mathbb{R}^n)$. A function $f:\mathbb{R}^n \rightarrow \mathbb{R}^m$ is Borel measurable if for each open set ${O} \subset \mathbb{R}^m$, the set $f^{-1}({O}) \coloneqq \{ x \in \mathbb{R}^n \; :\; f(x) \in {O}\}$ is Borel measurable \cite{rockafellar2009variational}. Further, we collect two useful results in martingale theory.
\begin{definition}[Supermartingale]\label{def:martingale}
A discrete-time supermartingale is a discrete-time stochastic process $X_0, X_1,X_2,\ldots$ that satisfies the following two properties for all $t$: (i) $\mathbb{E}[|X_t|] < \infty$, and (ii) $\mathbb{E}[X_{t+1}|X_1, \ldots, X_{t}] \leq X_t$.
\end{definition}
%Note that (ii) boils down to $\mathbb{E}[X_{t+1}| X_{t}] \leq X_t$ if $X_t$ evolves according to a Markovian dynamics.

\begin{lemma}[Ville's inequality, \cite{ville1939etude}]\label{lemma:Ville}
Let  $X_0, X_1,\ldots$ be a non-negative supermartingale. Then, for any  real $a > 0$,
$$
\emph{Pr}(\sup_{t \geq 0} X_t > a) \leq \frac{\mathbb{E}[X_0]}{a}.
$$
\end{lemma}

Let $\Sigma[x]$ denote the set of sum-of-squares polynomials, $\Sigma[x]^n$ the set of $n\times n$ sum-of-squares polynomial matrices, and $\mathbb{R}[x]^{n\times m}$ the set of $n\times m$ polynomial matrices with real coefficients. Finally, we recall a celebrated result in convex optimization, which is a direct consequence of the Positivstellens\"atze: \cite{parrilo2003semidefinite}.
 \begin{lemma}[S-Lemma]\label{eq:S-procedure1}
  Let \: $t_i(x)\in\Sigma[x], g_i(x)\in\mathbb{R}[x], f(x)\in\mathbb{R}[x], i \in \mathcal{I}$, then
$$
\begin{aligned}
& f(x)-\sum_{i \in \mathcal{I}}t_i(x)g_i(x)\in\Sigma[x]\Rightarrow \\
& f(x)\geq 0, \forall x\in\ \bigcap_{i \in \mathcal{I}}\{x:g_i(x)\ge 0\}.
\end{aligned}
$$
 \end{lemma}

%   \begin{lemma}[Farkas' Lemma]\label{eq:Farkas}
%   Let $A \in \mathbb{R}^{n \times m}$ and $b \in \mathbb{R}^m$. Then exactly one of the following assertions is true:
%   \begin{enumerate}
% \item There exists an $x \in \mathbb{R}^n$ such that $Ax = b$ and $x \geq 0$.
% \item There exists a $y \in \mathbb{R}^m$ such that $A^\top \!y \!\geq \!0$ and $b^\top\! y \! <\! 0$.
%   \end{enumerate}
%  \end{lemma}
 
 % Let $t(x)\in\Sigma[x]$, then
	% \begin{equation}\label{eq:S-procedure1}
	% 	f(x)-t(x)g(x)\in\Sigma[x]\Rightarrow f(x)\ge 0, \forall x\in\{x:g(x)\ge 0\}.
	% \end{equation}

%For a continuous function $b(x):\mathbb{R}^n \rightarrow \mathbb{R}$, a set denoted by the corresponding calligraphic letter $\mathcal{B}$ is defined by $\mathcal{B} \coloneqq \{x \in \mathbb{R}^n\; : \; b(x) \geq 0\}$.

\section{Problem setup}\label{sec:2}
Consider a controlled stochastic discrete-time system
\begin{equation}\label{eq:poly:sys}
x_{t+1} = Ax_t + Bu_t + Dw_t,
\end{equation}
with known (i.e., deterministic) initial state $x_0 \in \mathbb{R}^n$. Let $(\Omega, \mathcal{F}, {\mathbb{P}})$ be the probability space for the sequence $\mathbf{w}_\infty : \Omega \rightarrow \mathcal{W}^\infty$ of random variables, i.e., $\mathbf{w}_\infty = \{w_t\}_{t=0}^\infty$ for $w_t:\Omega \rightarrow \mathcal{W}$.  Let $\{\mathcal{F}_t\}_{t=0}^\infty$ denote the natural filtration for the sequence $\mathbf{w}_\infty$, that is $\mathcal{F}_t \subset \mathcal{F}$ consists of all sets of the form $\{\omega \in \Omega \: : \: \mathbf{w}_{[:t]} \in F\}$ for all $F \in \mathbf{B}(\mathcal{W}^t)$, where $\mathbf{w}_{[:t]}:\Omega \rightarrow \mathcal{W}^{t}, \mathbf{w}_{[:t]} = (w_0, w_1,\ldots, w_{t-1})$ is a subsequence. Informally speaking, $\mathcal{F}_t$ corresponds to the information available up to time $t$. For $1\le p \leq \infty$, let $\mathcal{L}_t^p$ denote the space of $\mathcal{F}_t$-measurable random variables $\mathbf{w}_{[:t]}$ with finite $p$-th moment, i.e., $\left(\mathbb{E}(|\mathbf{w}_{[:t]}|^p)\right)^{\frac{1}{p}} < +\infty$. We rely on the following assumption:

\begin{assumption}[Disturbance model]\label{ass:noise}
The disturbances $w_t : \Omega \rightarrow \mathcal{W}$ are independent and identically distributed (i.i.d.) in time, zero-mean sub-Gaussian random variables. Each random variable has an equivalent probability measure $\mu:\mathbf{B}(\mathcal{W}) \rightarrow [0,1]$ defined as $\mu(F) = \mathbb{P}(\{\omega \in \Omega : w_t(\omega) \in F\})$ for all $F \in \mathbf{B}(\mathcal{W}), \: t \geq 0$. The second
moment of $w_t$ is finite and the covariance matrix is defined as $\Sigma = \mathbb{E}[w_t w_t^\top]$ for all $t \geq 0$.
\end{assumption}

The dynamical system \eqref{eq:poly:sys} gives rise to a random process; we shall use $x(t,x_0,\mathbf{w}_{[:t]})$ to denote the random solution of \eqref{eq:poly:sys} at time $t$ given the deterministic initial condition $x_0$, i.e., $x_0(\omega) = x_0$ for all $\omega \in \Omega$, with $x_0 \in \mathcal{I}$.

 % We define a random process $\mathbf{x}$ as the sequence of random variables $x_t: \text{dom}(x_i) \subset \Omega \rightarrow \mathbb{R}^d$ and $\text{dom}(x_{t+1}) \subset \text{dom}(x_i)$. The random process $\mathbf{x}$ is \textit{adapted} to the natural filtration of $\mathbf{w}$ if $x_{t+1}$ is $\mathcal{F}_t$ measurable for each $t \geq 0$, i.e., if $x_{t+1}(\omega) = Ax_t(\omega) + B u(x_t(\omega)) + w_t(\omega)$ for all $\omega \in \text{dom}(x_{t+1})$ and $t \geq 0$. We use $x(t,x_0,\mathbf{w}_{[:t]})$ to denote the random solution of \eqref{eq:poly:sys} at time $t$ given the deterministic initial condition $x_0$, i.e., $x_0(\omega) = x_0$ for all $\omega \in \Omega$, with $x_0 \in \mathcal{I}$.

\subsection{Safety description}
We focus on safety requirements that are formulated as constraints on the system states as $x_t \in \mathcal{S}, \forall t \in \mathcal{T}$, where $\mathcal{T} = \{0, \ldots, T\}$ is the index set spanning the time horizon, and the \textit{safe} set $\mathcal{S}$ obeys the following definition:
\begin{definition}[Semi-algebraic set]
 A set $\mathcal{X} \subset \mathbb{R}^n$ is semi-algebraic if it can
be represented by polynomial equality and inequality constraints. If there are only equality constraints, the set is
called algebraic.
\end{definition}
Consequently, we consider a compact safe set defined as
\begin{equation}\label{eq:poly:set}
\mathcal{S}\coloneqq \bigcap_{j=1}^q \{x \in \mathbb{R}^n \; : \; s_j(x) \geq 0\},
\end{equation}
for given polynomials $s_1(x),\ldots, s_q(x)$. Semi-algebraic safe sets virtually encompass all safety specifications of engineering systems. We further consider a semi-algebraic initial set $\mathcal{I} \subseteq \mathcal{S}$ defined as
\begin{equation}\label{eq:poly:set:init}
\mathcal{I}\coloneqq \bigcap_{j=1}^p \{x \in \mathbb{R}^n \; : \; r_j(x) \geq 0\},
\end{equation}
for known polynomials $r_1(x),\ldots, r_p(x)$. Given the system \eqref{eq:poly:sys}, a safe set \eqref{eq:poly:set} and an initial set \eqref{eq:poly:set:init}, we call the system \textit{safe} if the closed-loop solution of \eqref{eq:poly:sys}, denoted as $x(t,x_0,\mathbf{w}_{[:t]})$,  satisfies either one of the following safety definitions:

\begin{definition}[Infinite-horizon safety]\label{def:bounded:safety}
Consider the system \eqref{eq:poly:sys}, a safe set $\mathcal{S}$, and a set of initial conditions $\mathcal{I}$. The system is safe if
\begin{equation*}
\bigcap_{j=1}^q s_j\left(x(t,x_0, \mathbf{w}_{[:t]})\right)\geq 0, \: a.s. \quad \forall t \in \mathcal{T}, T \rightarrow \infty
\end{equation*}
holds for all $x_0 \in \mathcal{I} \subseteq \mathcal{S}$.
\end{definition}
% \vspace{-0.5cm}
Definition \ref{def:bounded:safety} says that the system is safe if its state vector lies in the safe set almost surely for all time steps, i.e., for $T \rightarrow \infty$. Note that infinite-horizon safety can only be achieved if the support $\mathcal{W}$ is bounded. Conversely, if $\mathcal{W} = \mathbb{R}^d$, the system trajectory will leave the safe set with probability one in the long run. Clearly, this safety definition naturally extends to the case of nonstochastic (e.g., adversarial) disturbances with bounded magnitude.

\begin{definition}[$(1-\alpha)$-safety in probability]\label{def:finite:safety}
Consider the system \eqref{eq:poly:sys}, a safe set $\mathcal{S}$, a set of initial conditions $\mathcal{I}$. Fix $T \in \mathbb{N}$ finite. For $\alpha \in (0,1)$ the system is said to be $(1-\alpha)$-safe if it satisfies
$$
\emph{Pr}(\inf_{t \in \mathcal{T}} b(x_t) \geq 0) \geq 1 - \alpha$$
 for all $x_0 \in \mathcal{I} \subseteq \mathcal{S}$.
\end{definition}
% \vspace{-0.5cm}
Definition \ref{def:finite:safety} declares the system to be "safe" if the entire trajectory lies within the safe set with high enough probability. Note that this safety notion applies equivalently to the case where $\mathcal{W}$ is bounded. Let 
\begin{equation}\label{eq:finite:safety}
\begin{aligned}
&\text{Pr}_{\text{exit}}(T,x_0)\!=\!\!\text{Pr}\left(\!\inf_{t \in \mathcal{T}} \bigcap_{j=0}^q s_j\left({x}(t, x_0,\mathbf{w}_{[:t-1]})\right)  \!<\! 0 \!\right)\!.
\end{aligned}
\end{equation}
be the $T-$step exit probability, i.e., the probability that \eqref{eq:poly:sys} leaves the $0$-superlevel set of the safe set $\mathcal{S}$ within $T$ time-steps after starting at $x_0$. Then 
$$
\text{Pr}(\inf_{t \in \mathcal{T}} b(x_t) \geq 0) \geq 1 - \alpha \iff \text{Pr}_{\text{exit}}(T,x_0) \leq \alpha.
$$
% \begin{figure*}
% \centering
% \begin{minipage}[t]{0.45\textwidth}
%     \tikzset{every picture/.style={line width=0.75pt}} %set default line width to 0.75pt
%     \centering
%     \input{graphics/safety}
%     \caption{Example of safety notion. Left: in the bounded support case, system state trajectories can be ensured to lie in the safe indefinitely. Right: in the unbounded support case, a fraction $1-\alpha$ of $T-$long trajectories is certified to lie within the safe set for the entire task duration.}\label{fig:safety}
% \end{minipage}%
% \hfill
% \begin{minipage}[t]{0.45\textwidth}
%     \tikzset{every picture/.style={line width=0.75pt}} %set default line width to 0.75pt
%     \centering
%     \input{graphics/safety2}
%     \caption{Turning safety into a verification problem: If there exists a set $\mathcal{B}$ such that the forward invariance property and the set containment property are simultaneously satisfied, then the (deterministic) system is safe.}\label{fig:verif}
%     \end{minipage}
% \end{figure*}

\subsection{Problem formulation}
We will employ the framework of CBFs to certify the safety of stochastic discrete-time dynamical systems. Our problem statement is as follows:

\textbf{Problem.} For a given discrete-time stochastic system as in \eqref{eq:poly:sys}, safe set as in \eqref{eq:poly:set} and initial set as in \eqref{eq:poly:set:init}, provide a \underline{convex co-design} method for a CBF $b(\cdot)$ and a feedback controller $u(\cdot)$ that jointly guarantees safety of \eqref{eq:poly:sys} according to Def.~\eqref{def:bounded:safety} or \eqref{def:finite:safety}.

Following \cite{wang2024convex}, we will exploit the tight connection between safety and CBFs, as outlined next. Consider the system \eqref{eq:poly:sys}, a safe set $\mathcal{S}$, and an initial set $\mathcal{I \subseteq \mathcal{S}}$. We are looking for a set $\mathcal{B} = \{ x \in \mathbb{R}^n \; :\; b(x) \geq 0$\}, where $b: \mathbb{R}^n \rightarrow \mathbb{R}$ is a continuous function, satisfying (i) a \textit{forward invariance property}, i.e., if $b(x_t) \geq 0$ then $b(x_{t+1}) \geq 0$ in an "appropriate" sense (i.e., with respect to the stochasticity handling); and (ii) a \textit{set containment condition}, i.e., $ \mathcal{I} \subseteq \{ x \in \mathbb{R}^n \; : \; b(x) \geq 0 \} \subseteq \mathcal{S}$. If both conditions hold simultaneously, we refer to $b(\cdot)$ as a \textit{control
 barrier function}, and to $u(\cdot)$ as a \textit{safe controller}. As a result, we can interpret the safety enforcement problem into an \textit{existence} problem: if there exists a set $\mathcal{B}$ (which in general depends on the controller $u(x)$) satisfying the above properties, then the system is safe.

We comment on the non-triviality of our problem statement with the following pedagogical example. Consider the simplified case of a deterministic system, e.g. $D = 0$ in \eqref{eq:poly:sys}. The forward invariance condition reads $ b(Ax + Bu(x)) - b(x) \geq -\beta(b(x)), \: \forall x \in \mathcal{B}.$ Since $b(\cdot)$ and $u(\cdot)$ are both parametrized functions, the composition $b \circ (Ax + Bu(\cdot))$ is in general nonconvex in the parameters that parametrizes $b$ and $u$. We highlight that this nonconvexity does not arise in the continuous-time regime for the popular case of input-affine systems, where the forward invariance condition is instead $\frac{\partial b(x)}{\partial x}f(x,u) \geq 0, \: \forall \: x \in \partial \mathcal{S}$ leading to a bilinearity that is easier to handle \cite{wang2023assessing}. As explained in Section 1, the majority of the literature considers the CBF to be given \emph{a priori}, often assuming\footnote{This can lead to pathological behaviors \cite[Case~1]{wang2024convex}.} $b(\cdot) = \min_{j = 1,\ldots,q} s_j(\cdot)$. A different option is to parameterise the CBF as a first-degree polynomial $b(x) := \ell x + h$ such that the composition boils down to $\ell \cdot b(Ax + u(x)) +  h$ \cite{wang2023assessing}. The resulting bilinearity between $\ell$ and $u(\cdot)$ can at this point be addressed by alternating between two sum-of-squares programs \cite{tan2004searching}, inheriting the drawbacks outlined in Section 1. Further, the predefined choice of a first-order polynomial might be conservative (and even lead to infeasibility).

\section{Co-design of CBFs and state feedback controllers}\label{sec:main}
In this section, we propose convex formulations to co-design a CBF $b(x)$ and a safe feedback controller $u(x)$ for system \eqref{eq:poly:sys} for different geometric shapes of the initial set $\mathcal{I}$ and safe set $\mathcal{S}$. We will consider the following parameterisation:
\begin{subequations}\label{eq:parametrization:all}
\begin{align}
b(x) & = 1 - x^\top \Omega^{-1}x,\\
u(x) & = Y \Omega^{-1}x
\end{align}
\end{subequations}
where $\Omega \in \mathbb{R}^{n \times n}, Y \in \mathbb{R}^{m \times n}$ are matrices to be designed. Notice that higher-order CBFs would potentially yield a larger invariant set $\mathcal{B}$ (hence reducing conservatism); however, they suffer from two drawbacks. First, they return a nonlinear design problem, which is not obviously convexifiable in general. Second, they require knowledge of higher-order moments of the probability distribution, which is often limiting in real-world applications (e.g., when moments are estimated from historical data). As a result, we proceed with a quadratic parametrization for $b(x)$ as it offers a good compromise between expressivity and computational tractability.
To simplify presentation, we derive our main result based on a specific geometry of $\mathcal{S}$ and $\mathcal{I}$ and later extend it.
\vspace{-0.2cm}
\begin{assumption}[Initial and safe set]\label{ass:sets}
$\mathcal{S}$ is non-empty, and
 $s_j(x)=a_j^\top x+1,j=1,\ldots,q.$ Further, $\mathcal{I} = \{ x\in\mathbb{R}^n: 1 - x^\top R x \geq 0\}$ for a given $R = R^\top \succ 0$.
\end{assumption}
We defer all the proofs of the subsequent results to the Appendix.

\subsection{Design for infinite-horizon safety}\label{section:infinite}
We consider uncertain systems as in \eqref{eq:poly:set} with $\mathcal{W}$ bounded and obeying the following structure, which we will generalize in Section 4:

\begin{assumption}[Bounded support]\label{eq:support}
The disturbance is bounded in $\mathcal{W}:= \{ w \in \mathbb{R}^d : w^\top w \leq 1\}$.
\end{assumption}

When confronted with disturbances, we consider a robust version of forward invariance reading
\begin{equation}\label{eq:robust:c1}
(w^\top w \leq 1) \wedge (x \in \mathcal{B}) \implies (Ax + Bu(x) + Dw) \in \mathcal{B}.
\end{equation}
A \textit{sufficient} condition for \eqref{eq:robust:c1} to hold is provided by the following \textit{global} implication: For all $x \in \mathbb{R}^n$ and a fixed $\beta \in (0,1)$, it holds
\vspace{-0.2cm}
\begin{equation}\label{eq:robust:c2}
(w^\top w \leq 1) \implies b(Ax + Bu(x) + Dw) - b(x) \geq -\beta b(x).
\end{equation}
Despite providing a conservative approximation of \eqref{eq:robust:c1}, \eqref{eq:robust:c2} allows for a convex reformulation, as shown next. Fix $\lambda, \beta \in (0,1)$ and let $\tilde{\beta} = 1 - \beta \in (0,1)$. Further, consider the following \textit{jointly convex} optimization problem:
\begin{subequations}\label{eq:convex-without-input-constraints}
\begin{align}
    \max~&\mathrm{log det}(\Omega)\\
    \mathrm{s.t.}~&  Y \in\mathbb{R}^{m\times n},  0\prec\Omega=\Omega^\top\in\mathbb{R}^{n\times n}, \\ 
    &\begin{bmatrix}
        -\tilde{\beta} \Omega + \lambda \Omega & 0 & \Omega A^\top + Y^\top B^\top \\
        0 & -\lambda I & D^\top\\
        A\Omega + BY & D & -\Omega
    \end{bmatrix} \preceq 0,\label{eq:convex-e}\\
     &\begin{bmatrix}
     R&I\\I&\Omega
     \end{bmatrix}\succeq 0,\label{eq:convex-f}\\
    &1-a_j^\top \Omega a_j\geq 0,\: j=1,\ldots,q,\label{eq:convex-h}
\end{align}
\end{subequations}
where $\Omega$, $Y$ are decision variables. We show next that \eqref{eq:convex-without-input-constraints} provides a convex design procedure for the synthesis of a valid CBF and feedback controller for \eqref{eq:poly:sys}, certifying infinite-horizon safety according to Definition \eqref{def:bounded:safety}. 

% \begin{tcolorbox}[colback=purple!6, colframe=white, sharp corners=south, boxrule=0pt,boxsep=-0.5mm]
\begin{theorem}\label{th:convex-design}
    Consider the system \eqref{eq:poly:sys}, safe set $\mathcal{S}$, and initial set $\mathcal{I}$. Let Assumptions \ref{ass:noise} and \ref{ass:sets} hold, and let $\mathcal{U}=\mathbb{R}^m$. Assume that for a fixed $\lambda \in (0,1)$ a solution to \eqref{eq:convex-without-input-constraints} exists and is denoted by $\Omega,Y$. Set $u = u(x)=Y\Omega^{-1}x$ and define $\mathcal{B} = \{ x \in \mathbb{R}^n \; : \; 1 - x^\top \Omega^{-1}x \geq 0\}$. Then the following hold
    \begin{enumerate}
        \item $\mathcal{B}$ is an invariant set for $x_{t+1} =Ax_t+Bu(x_t) + Dw$.
        \item $\mathcal{I}\subseteq \mathcal{B}\subseteq\mathcal{S}$.
        \item $\mathcal{B}$ is the set with largest volume satisfying 1) and 2) given the parameterisation in \eqref{eq:parametrization:all}.
    \end{enumerate}
    Further, system \eqref{eq:poly:sys} under the controller $u(x)$ is safe according to Definition \eqref{def:bounded:safety}.
\end{theorem}
% \end{tcolorbox}
% Note that the objective function in \eqref{eq:convex-without-input-constraints} can be substituted by $\max \mathrm{Tr}(\Omega)$ \cite{durieu1996trace}, that is better supported by available SOS parsers \cite{prajna2002introducing}.

\begin{remark}[On optimizing over $\lambda$ and $\beta$]
Our procedure assumes a fixed $\lambda \in (0,1]$. If the control designer wants to optimize over $\lambda$ as well, constraint \eqref{eq:convex-e} becomes nonconvex due to the existence of the bilinear matrix term $\lambda \Omega$. Let $g(\lambda)$ denote the optimal value of \eqref{eq:convex-without-input-constraints} for a fixed value of $\beta$ (and hence $\tilde{\beta}$). We immediately notice that $g(\lambda)$ is not well-defined for $\lambda\leq 0$ or $\lambda > 1$ as the optimization program becomes infeasible (see also proof of Theorem \eqref{th:convex-design}). When $\lambda \in (0,1]$, $g(\lambda)$ can be shown to be quasi-convex. Thus, the optimization over $\lambda$ can be efficiently carried out by bisection on $(0,1]$ to find $\max_\lambda g(\lambda)$.  Similar reasoning can be applied to $\beta \in (0,1]$.
\end{remark}

 Finally, note that in the case of an ellipsoidal safe set of the form $\mathcal{S} = \{ x \in \mathbb{R}^n \; \; 1 - x^\top S x \geq 0\}$ with $S \succ 0$ given, the containment condition $\mathcal{B} \subseteq \mathcal{S}$ in \eqref{eq:convex-h} amounts to
$$
\begin{aligned}
& \{x \in \mathbb{R}^n \:: \: 1 - x^\top \Omega^{-1}x \geq 0\} \subseteq \{x \in \mathbb{R}^n \:: \: 1 - x^\top S x \geq 0\}\\
& \hspace{2.7cm}\iff  \begin{bmatrix} \Omega & \Omega \\ \Omega &  S^{-1} \end{bmatrix} \succeq 0. 
\end{aligned}
$$

\subsection{Design for finite-horizon safety}\label{sec:design:finite}
We next consider the case of uncertainty with unbounded support $\mathcal{W}$. In this case we relax our notion of safety to the requirement:
\begin{equation}\label{eq:safety:finite:recap}
\begin{aligned}
&\text{Pr}(x(t, x_0, \mathbf{w}_{[:t-1]}) \in \mathcal{S}, \forall t = \{0,\ldots,T\})\\
= & 1 - \text{Pr}\left(\min_{t \in \{0,\ldots,T\}} b_t(x(t, x_0, \mathbf{w}_{[:t-1]})) < 0 \right)\\
\geq & 1 - \alpha,
\end{aligned}
\end{equation}
for $\alpha \in (0,1)$ being a user-defined risk tolerance parameter. We will exploit tools from martingale theory to encode \eqref{eq:safety:finite:recap} via a convex optimization problem. 

Martingale theory provides sufficient conditions for \eqref{eq:safety:finite:recap} to hold in the form of constraints on the expected value of $b(\cdot)$ of the form
\begin{equation}\label{eq:safety:martingale}
\mathbb{E}[b(Ax_t + Bu(x_t) + w_t) | \mathcal{F}_{t}]  \geq \tilde{\beta} b(x_t) + \delta,
\end{equation}
with $\tilde{\beta} = 1 - \beta, \: \beta \in (0,1)$ for all $t \in \mathcal{T}$. Note that in \eqref{eq:safety:martingale} we have introduced the additional parameter $\delta$ for greater generality: in particular, $\delta < 0$ would reduce the conservatism of our design, while $\delta > 0$ further constrains it. We begin by considering a fixed $\delta \in \mathbb{R}$ and later show how to harness it to impose precise safety requirements. Before providing our main result of this subsection, we describe the proposed procedure consisting of three main steps. First, by recognizing that the sequence $\{b_t\}$ with $b_t \overset{\triangle}{=} b(x_t)$ might not have the desired martingale properties, we define a non-negative supermartingale $\zeta_t$ starting from the increase condition in \eqref{eq:safety:martingale} by shifting and scaling it appropriately. Second, by noticing that Ville's inequality (Lemma \eqref{lemma:Ville}) is now applicable to the constructed $\zeta_t$, we define a parameter $a^\star \in \mathbb{R}$ such that
$$
\left(\sup_{0 \leq t \leq T} \zeta_t \leq  a^\star \right) \implies \left( \inf_{0 \leq t \leq T} b(x_t) \geq 0 \right)
$$
in order to map the probabilistic guarantees back to the original system (e.g., in $b(\cdot)$). Third, armed with the bound derived from the previous step, we embed the condition in \eqref{eq:safety:martingale} into a convex optimization problem, while restricting the design to match with the user-defined maximum tolerance $\alpha$.

We consider each step individually. The following lemma details criteria to  construct a non-negative supermartingale starting from \eqref{eq:safety:martingale}.
\begin{lemma}\label{lemma:super:general}
Fix $\beta \in (0,1), \delta \in (\beta - 1, \beta]$ and define $\psi = \beta - \delta$. Let \eqref{eq:safety:martingale} hold for all $t  \in \mathcal{T}$ with $T \in \mathbb{N}$ finite. Consider functions $g_i(\beta, \delta)$ and $h_i(\beta, \delta)$ such that i) $0 < g_i(\beta, \delta) \leq \frac{1}{1-\beta}$ for all $i$, ii) $h_i(\beta,\delta) \geq 0$ for all $i$, and iii) $\prod_{i=1}^{t+1} g_i(\beta, \delta) \psi \leq h_{t+1}(\beta, \delta) - h_{t}(\beta, \delta)$ for all $t$. Then 
$$
\zeta_t \coloneqq \left( \prod_{i=1}^t g_i(\beta,\delta)\right)x_t^\top \Omega^{-1}x_t + \left( H - \sum_{i=1}^t h_i(\beta,\delta) \right)
$$
is a non-negative supermartingale in the interval $t\in\mathcal{T}$ for any $H \geq \sum_{i=1}^T h_i(\beta, \delta)$. 
\end{lemma}
\vspace{-0.2cm}
Intuitively, $g_i$ cancels out the effect of the multiplicative term $\tilde{\beta}$ and $h_i$ the effect of $\delta$. Next, we show how the nonnegative supermartingale defined in Lemma \eqref{lemma:super:general} can be used together with Lemma \eqref{lemma:Ville} to ensure the sought probabilistic safety guarantees.
\vspace{-0.2cm}
\begin{proposition}\label{prop:Ville:general}
Fix $\beta \in (0,1), \delta \in (\beta-1, \beta]$ and define $\psi = \beta - \delta$. Let \eqref{eq:safety:martingale} hold for all $t  \in \{0,\ldots, T\}$. Further, let $\mathcal{B} = \{ x \in \mathbb{R}^n \; : \; 1 - x^\top \Omega^{-1}x \geq 0\}$ and consider any functions $g_i(\beta, \delta)$ and $h_i(\beta, \delta)$ satisfying Lemma \ref{lemma:super:general}. Define $a(t) = \left( \prod_{i=1}^t g_i(\beta, \delta) \right) + \left( H - \sum_{i=1}^t h_i(\beta, \delta) \right)$ and let $a^\star = \inf_{t \in \mathcal{T}} a(t)$. Then
$$
\emph{Pr}(x_t \in \mathcal{B}, \forall t \in \mathcal{T}) \geq 1 - \frac{\mathbb{E}[\zeta_0]}{a^\star}.
$$
\end{proposition}
\vspace{-0.2cm}
For any suitable choice of functions $g_i(\beta,\delta), h_i(\beta,\delta)$, Proposition \ref{prop:Ville:general} quantifies the joint-in-time safety probability of system \eqref{eq:poly:sys}. Consequently, the control engineer might want to design them to achieve the highest possible bound on the safety probability. In the next example, we propose a constructive design of these elements, inspired by \cite{kushner1967stochastic,cosner2023robust}.

\begin{example}[Constructive design of $\zeta_t$]
Consider the setting of Lemma \ref{lemma:super:general}. For all $i = \{0,\ldots,T\}$, we let $g_i(\beta,\delta) = \eta$ with $1 < \eta \leq \frac{1}{1-\beta}$ and  $h_i(\beta,\delta) = \psi \eta^i$. We further define $H = \psi \sum_{i=0}^T \eta^i$. Notice that properties i) and ii) are trivially satisfied given the bounds on $\eta$, while property iii) holds since
$$
\eta^{t+1} \psi + \psi \left( \sum_{i=1}^{T}\eta^i - \sum_{i=1}^{t+1} \eta^i\right) = \psi \left( \sum_{i=1}^{T}\eta^i - \sum_{i=1}^{t} \eta^i\right).
$$
Then, starting from
$
\mathbb{E}[x_{t+1}^\top \Omega^{-1}x_{t+1}] \leq \tilde{\beta} (x_t^\top \Omega^{-1} x_t) + \psi,
$
and following similar steps as in the proof of Lemma \ref{lemma:super:general} we get
\begingroup
\allowdisplaybreaks
$$
\begin{aligned}
& \mathbb{E}\left[ \eta^{t+1}(x_{t+1}^\top \Omega^{-1}x_{t+1}) + \psi\left(\sum_{i=1}^T \eta^{i} - \sum_{i=1}^{t+1} \eta^{i} \right)\right] \\
\leq & \eta^{t+1} \left(\tilde{\beta} (x_t^\top \Omega^{-1} z_t) + \psi\right)+ \psi \left(\sum_{i=1}^T \eta^{i} - \sum_{i=1}^{t+1} \eta^{i} \right)\\  
\aleq & \eta^t (x_t^\top \Omega^{-1} x_t) +  \eta^{t+1}\psi + \psi \sum_{i=1}^T \eta^{i} - \psi\sum_{i=1}^{t+1} \eta^{i}\\
\beq & \eta^t (x_t^\top \Omega^{-1} x_t) + \psi \sum_{i=1}^T \eta^{i} - \psi\sum_{i=1}^{t} \eta^{i}\\
\ceq & \eta^t (x_t^\top \Omega^{-1} x_t) + \psi \eta \frac{\eta^T - \eta^t}{\eta - 1}=: \zeta_t,
\end{aligned}
$$
\endgroup
where $(a)$ follows by property i), $(b)$ by property iii) and $(c)$ from the relationship $\eta \frac{\eta^T - \eta^t}{\eta - 1} = \sum_{i=1}^T \eta^{i-1} - \sum_{i=1}^t \eta^{i}$. It is easy to see that $\zeta_t$ is a non-negative supermartingale since $\eta > 1, \psi \geq  0 $ by assumption, and $\Omega^{-1} \succ 0$. By definition of $a(t)$ in Proposition \ref{prop:Ville:general}, we have
\begin{equation}\label{eq:a:def}
a(t) = \eta^t ( \underbrace{1 - \frac{\psi \eta}{\eta - 1}}_{\kappa})+ \frac{\psi \eta}{\eta - 1}\eta^T.
\end{equation}
We consider two cases. First, let $\delta < 0$. In this case $\psi = \beta - \delta > \beta$, hence $\frac{1}{1-\beta} < \frac{1}{1-\psi}$ and consequently $\eta \leq \frac{1}{1-\psi}$ for any $\eta \in (1, \frac{1}{1-\beta}]$.
As a result, $\kappa < 0$ for $\eta \in (1, \frac{1}{1-\beta}]$, implying that $a(t)$ is monotonically decreasing, corresponding to the case where $\frac{\partial a(t)}{\partial t}<0$. Its minimum is achieved at $t = T$, yielding
\begin{equation}\label{eq:first}
\alpha_1(\eta) = \frac{1-b(x_0)+ \frac{\psi \eta}{\eta-1}(\eta^T - 1)}{\eta^T}.
\end{equation}
Note that $\alpha_1(\eta)$ is a decreasing function of $\eta$ since $\frac{\partial \alpha_1(\eta)}{\partial \eta} 
\leq 0$ and hence the tightest bound is achieved by choosing $\eta^\star = \frac{1}{1-\beta}$.
Let now $\delta \geq 0$. In this case, $\psi \leq \beta$ and consequently $\frac{1}{1-\psi} \leq \frac{1}{1-\beta}$. In the interval $\eta \in [ \frac{1}{1-\psi} \leq \frac{1}{1-\beta}]$, it holds $\kappa \geq 0$, corresponding to the case where $\frac{\partial a(t)}{\partial t}>0$, and thus $a^\star$ is achieved at $t = 0$, yielding
\begin{equation}\label{eq:second}
\alpha_2(\eta) = 1 - \frac{b(x_0)}{1 + \frac{\psi \eta}{\eta - 1}(\eta^T - 1)}.
\end{equation}
The tightest bound is in this case achieved by $\eta^\star = \frac{1}{1-\psi}$ since $\frac{\partial \alpha_2(\eta)}{\partial \eta} 
\geq 0$. In the interval $\eta \in (1, \frac{1}{1-\psi})$ we fall back to \eqref{eq:first}.
 $\hfill \square$
\end{example}

\vspace{-0.2cm}
We are now left with the task of encoding \eqref{eq:safety:martingale} in a convex optimization program. To simplify the presentation, we begin by assuming that both $\beta$ and $\delta$ are fixed. Further, we restrict our attention to the design rules in Example 1; similar conclusions can be derived for different design choices obeying Lemma \ref{lemma:super:general}. Fix $\beta \in (0,1), \delta \in (\beta-1,\beta]$ and let $\tilde{\beta} = 1 - \beta$. In line with Assumption \ref{ass:sets}, consider an ellipsoidal initial set $\mathcal{I} = \{x \in \mathbb{R}^n \; : \; 1 - x^\top R x \geq \sigma\}$, where $R = R^\top \succ 0$ and $\sigma \in (0,1)$ are given. Further, consider the following \textit{jointly convex} optimization problem:
\begin{subequations}\label{eq:convex-without-input-constraints-finite}
\begin{align}
    \max~&\mathrm{logdet}(\Omega)\\
    \mathrm{subject~to}~& \lambda \geq 0,\\
    &0\prec\Omega=\Omega^\top\in\mathbb{R}^n, Y\in\mathbb{R}^{m\times n} \label{eq:convex-fin-b}\\
    &\begin{bmatrix}
        \lambda I& \Sigma^{1/2}\\
        \star&\Omega
    \end{bmatrix}\succeq 0,\label{eq:convex-fin-n}\\ 
    & \beta - \delta - \lambda \geq 0,\label{eq:convex-fin-m}\\
    &\begin{bmatrix}
        \tilde{\beta}\Omega&\Omega A^\top +B^\top Y^\top \\
        \star&\Omega
    \end{bmatrix}\succeq 0,\label{eq:convex-fin-e}\\
     &\begin{bmatrix}
     \frac{R}{1-\sigma}&I\\I&\Omega
     \end{bmatrix}\succeq 0,\label{eq:convex-fin-f}\\
   % & 1 + \theta_1 \delta + \theta_2 \mu \leq \alpha,\\
    &1-a_j^\top \Omega a_j\ge 0,j=1,\ldots,q,\label{eq:convex-fin-h}
\end{align}
\end{subequations}
where $\Omega, Y, \lambda$ are decision variables. The next theorem provides a bound on the safety probability of the trajectories generated by system \eqref{eq:poly:sys} starting from $\mathcal{I}$.

%s, while $\theta_1 = \frac{(1-\beta)^{-T} - 1}{\beta(1-\beta)^{-T}}, \theta_2 = \frac{1}{(1-\beta)^{-T}}$ are known constant. The following theorem constitutes the main result of this subsection. %We show next that \eqref{eq:convex-without-input-constraints} provides a convex design procedure for the synthesis of a valid CBF and feedback controller for \eqref{eq:poly:sys}, certifying infinite-horizon safety according to Definition \eqref{def:bounded:safety}. 

%\begin{tcolorbox}[colback=purple!6, colframe=white, sharp corners=south, boxrule=0pt,boxsep=-0.5mm]
\begin{theorem}\label{th:convex-design:pp}
    Consider the system \eqref{eq:poly:sys}, safe set $\mathcal{S}$, and initial set $\mathcal{I}$. Let Assumptions \eqref{ass:noise} and $\mathcal{U}=\mathbb{R}^m$. Assume that a solution to \eqref{eq:convex-without-input-constraints-finite} exists. Set $u = u(x)=Y\Omega^{-1}x$ and let $\psi = \beta -\delta$ where $\beta\in(0,1),\delta\in(\beta-1,\beta].$ Then the system \eqref{eq:poly:sys} under $u(x)$ is $(1-\alpha)$-safe according to Definition \eqref{def:finite:safety} for all $x_0 \in \mathcal{I}$, with
    $$
    \alpha \leq \begin{cases}
\frac{1-\sigma+ \frac{\psi \eta}{\eta-1}(\eta^T - 1)}{\eta^T} & \text\:{if}\: \delta < 0,\\
1 - \frac{\sigma}{1 + \frac{\psi \eta}{\eta - 1}(\eta^T - 1)} & \text{otherwise.}
    \end{cases}
    $$
\end{theorem}
%\end{tcolorbox}

While Theorem \ref{th:convex-design:pp} provides a tool to upper-bound the exit probability, it does not allow us to constrain it below a certain user-specified threshold as the parameters influencing $\alpha$ are \textit{given}. In the following, we show how to exploit $\delta$ as a degree of freedom to impose conditions ensuring a prescribed user-defined $\bar{\alpha} \in (0,1)$. To this end, recall the design from Example 1. Fix $\beta \in (0,1)$; in turn, $\delta$ is allowed to lie in the interval $(\beta-1,\beta]$. Then by evaluating \eqref{eq:first} and \eqref{eq:second} at their respective $\eta^\star$ to achieve the highest bound possible, we get 
$$
\begin{aligned}
& \text{Pr}_{\text{exit}}(T, x_0)\\
\leq & \begin{cases}
(1 - \sigma)(1-\beta)^T + (\beta - \delta) \sum_{i=1}^T (1-\beta)^{i-1} & \text{if}\: \delta < 0\\
1 - \sigma(1-\beta+\delta)^T & \text{oth.}
\end{cases}
\end{aligned}
$$
We require $\text{Pr}_{\text{exit}}(T, x_0) \leq \bar{\alpha}$. Defining the constants $\theta_1 := (1-\beta)^T, \theta_2 := (1-\beta)^T + \beta\sum_{i=1}^T (1-\beta)^{i-1}, \theta_3 := \sum_{i=1}^T (1-\beta)^{i-1}$ and $\theta_4:=\sqrt[T]{\frac{1 - \alpha}{\sigma}}$, this relationship can be encoded as a set of affine conditions in $\delta$:
$$
\begin{aligned}
\text{Pr}_{\text{exit}}(T, x_0) \leq & \begin{cases}
 \theta_3 \delta \geq \theta_2 - \theta_1 \sigma - \alpha & \text{if}\: \delta < 0\\
\delta \geq \theta_4 + \beta - 1  & \text{oth.}
\end{cases}
\end{aligned}
$$
We encode the above cases via a big-M reformulation making use of an auxiliary binary variable $z \in \{0,1\}$, where $z = 1$ if $\delta \geq 0$ and $z = 0$ if $\delta < 0$. Let $M > 0$ be a "large enough" constant and consider the following set of conditions: 
\begin{equation}\label{eq:bigM}
\begin{cases}
\delta \in (\beta-1,\beta], z \in \{0,1\}\\
\delta \geq -M(1-z)\\
\delta \leq M z,\\
\theta_2 - \theta_1 \sigma - \alpha - \theta_3 \delta \leq M(1-z),\\
\theta_4 + \beta - 1 - \delta \leq Mz.
\end{cases}
\end{equation}
Let $\Delta = \{ \delta \: | \: \eqref{eq:bigM} \: \text{is satisfiable for some $z$}\}$. Following a similar reasoning as before, it is easy to see that adding the constraint $\delta \in \Delta$ to Problem \ref{eq:convex-without-input-constraints-finite} yields $\text{Pr}(x_t \in \mathcal{S}, \: \forall t \in \mathcal{T}) \geq 1 - \bar{\alpha}$. At this point, \ref{eq:convex-without-input-constraints-finite} amounts to solving two SDPs, one with $z = 0$ and one with $z = 1$. 

\section{Extensions}\label{sec:4}
In this section we provide several extensions to the framework presented in Section \ref{sec:main}, greatly expanding the applicability of our results.

\subsection{Dealing with input constraints}
In Section \ref{sec:main} we showed how to co-design a CBF and a feedback controller, assuming that $\mathcal{U}=\mathbb{R}^m$. We now extend the results to the case that the control ability is constrained, by providing convex conditions ensuring $u=u(x)$ to be bounded. For simplicity, consider the bounded support setting of Subsection \ref{section:infinite}. We require the following implication to hold
\begin{equation}\label{eq:input:implication}
x \in \mathcal{B} \implies u(x) \in  \mathcal{U}.
\end{equation}
We provide convex sufficient conditions for \eqref{eq:input:implication} to hold for different types of the input set $\mathcal{U}$.

\begin{lemma}\label{lemma:input:constraint}
    Consider system \eqref{eq:poly:sys}, $\mathcal{B}:=\{x\in\mathbb{R}^n:1-x^\top \Omega^{-1}x\ge 0\}$, and let $\mathcal{U}_1:=\{u\in\mathbb{R}^m:Hu\le h\}$, where $\mathcal{H}\in\mathbb{R}^{k\times n}$, $h\in\mathbb{R}^{k}$. If there exists $\lambda>0$, and $Y\in\mathbb{R}^{m \times n}$, such that
    \begin{equation}\label{eq:convex-lin-input}
        \begin{bmatrix}
            \mathcal{H}_{11}^i&\mathcal{H}_{12}\\\mathcal{H}_{12}^\top &I_{n+1}
        \end{bmatrix} \succeq 0,\: i=1,\ldots,k,
    \end{equation}
    where
    \begin{equation*}
        \mathcal{H}_{11}^i=\begin{bmatrix}
            \Omega&Y^\top H_i^\top \\H_iY&2\lambda h_i
        \end{bmatrix},
        \mathcal{H}_{12}=\begin{bmatrix}
            0&0\\0&\lambda
        \end{bmatrix}.
    \end{equation*}
    Then, $u=u(x)=Y\Omega^{-1}x\in\mathcal{U}_1$, for any $x\in\mathcal{B}$.
\end{lemma}
% \begin{proof}
%  Equation \eqref{eq:convex-lin-input} is sufficient to
%     \begin{equation*}
%         \mathcal{H}_{11}^i-\mathcal{H}_{12}I_{n+1}\mathcal{H}_{12}^\top \succeq 0,i=1,\ldots,k,\forall x\in\mathbb{R}^n,
%     \end{equation*}
%     which implies
%     \begin{equation*}
%         \Omega-\frac{Y^\top H_i^\top H_iY }{2\mu h_i-\mu^2}\succeq 0,i=1,\ldots,k,\forall x\in\mathbb{R}^n.
%     \end{equation*}
%     Given that $\mu >0$, we have
%     \begin{equation*}
%         \mu\Omega-\frac{Y^\top H_i^\top H_iY}{2h_i-\mu}\succeq 0,i=1,\ldots,k,\forall x\in\mathbb{R}^n.
%     \end{equation*}
%     Left- and right- multiply the matrix by $\Omega^{-1}$, and define $K=Y\Omega^{-1}$, we obtain
%     \begin{equation*}
%         \mu\Omega^{-1}-\frac{K^\top H_i^\top H_iK}{2h_i-\mu}\succeq 0,i=1,\ldots,k,\forall x\in\mathbb{R}^n.
%     \end{equation*}
%     Apply Schur complement again, we have for all $x\in\mathbb{R}^n$
%     \begin{align*}
%         &\begin{bmatrix}
%             \mu\Omega^{-1}&-K^\top H_i^\top \\-H_iK&2h_i-\mu
%         \end{bmatrix}\succeq 0,i=1,\ldots,k,\forall x\in\mathbb{R}^n \Longleftrightarrow\\
%        &-2(H_iKx-h_i)+\mu(x^\top \Omega^{-1}x-1)\ge 0,i=1,\ldots,k..
%     \end{align*}
%     Then we have for $i=1,\ldots,k$, $H_iu(x)\le h_i$ for any $x$ such that $b(x)\ge 0$. Hence, we conclude the proof. 
% \end{proof}

Note that \eqref{eq:convex-lin-input} is convex in $\Omega$, $Y(x)$, and $\mu$. Therefore, it can be incorporated into \eqref{eq:convex-without-input-constraints} retaining convexity of the resulting design procedure. As a special case of Lemma \ref{lemma:input:constraint}, one can enforce input constraints of the form $u \in \mathcal{U}_2$, with $\mathcal{U}_2 \{ u \in \mathbb{R}^m \; : \; \|u\|_\infty \leq \bar{u}\}$ above rewrite the infinity-norm as a series of affine inequality constraints. Finally, we can enforce 2-norm constraints, following \cite{wang2024convex}.
\begin{lemma}
    Consider system \eqref{eq:poly:sys}, and let $\mathcal{U}_2:=\{u\in\mathbb{R}^m:u^\top u\le \bar{u}\}$, where $\bar{u} > 0$. If there exists $Y \in \mathbb{R}^{m \times n}$, such that
    \begin{equation}
        \begin{bmatrix}
            \Omega&Y^\top \\Y&\bar{u}
        \end{bmatrix}\succeq 0.
    \end{equation}
     Then, $u=u(x)=Y\Omega^{-1}x\in\mathcal{U}_2$, for any $x\in\mathcal{B}$.
\end{lemma}

We conclude by discussing the implications of these results for the setting in Subsection \ref{sec:design:finite}. Since $\text{Pr}(x_t \in \mathcal{B}, \forall t \in \mathcal{T}) \geq 1 - \alpha$ for $T \in \mathbb{N}$ finite, by enforcing \eqref{eq:input:implication} we correspondingly get $\text{Pr}(u_t \in \mathcal{U}, \forall t \in \mathcal{T}) \geq 1 - \alpha$.

\subsection{Dealing with distributional ambiguity}\label{sec:ambiguity}
The design procedure outlined in Subsection \ref{sec:design:finite} assumes exact knowledge of the true probability measure $\mu$ (compare (14d)). However, in practice, one might not have such information available. On the contrary, the control designer might only know a \textit{nominal} measure $\hat{\mu}$ which inevitably differs from the true one. Therefore, it is essential to robustify the presented design procedure against this mismatch. Let $\hat{\mu}$ and $\hat{\mathcal{W}}$ be estimates of the true probability measure $\mu$ and support $\mathcal{W}$. We consider two sources of ambiguity: (i) unmodeled disturbances $(\hat{\mathcal{W}} \neq \mathcal{W})$ with $\hat {\mathcal{W}}$ bounded. and (ii) incorrectly modeled disturbances $(\hat{\mu} \neq \mu)$. Notice that (i) influences the infinite-time safety as per Def.~\ref{def:bounded:safety} , while (ii) influences the $(1-\alpha)$ safety as per Def.~\ref{def:finite:safety}.

We begin with the first case. Let $\hat{\mathcal{W}} = \{ w \in \mathbb{R}^d \; : \; w^\top w \leq 1\}$ according to Assumption \ref{eq:support} and assume without loss of generality that $\hat{\mathcal{W}} \subseteq {\mathcal{W}}$. The Hausdorff distance between the two sets evaluates
\begin{equation}
\begin{aligned}
d_H({\mathcal{W}}, \hat{\mathcal{W}}) & := \max \{ \sup_{w \in \mathcal{W}} |x|_{\hat{\mathcal{W}}}, \sup_{w \in \hat{\mathcal{W}}} |x|_{{\mathcal{W}}} \}\\
& = \max_{w \in \mathcal{W}} |x|_{\hat{\mathcal{W}}}.
\end{aligned}
\end{equation}
Let $d \overset{\triangle}{=} d_H({\mathcal{W}}, \hat{\mathcal{W}})$ to shorten the notation. We can then inflate $\tilde{\mathcal{W}} = \hat{\mathcal{W}} + d \mathcal{E}$ where $\mathcal{E}$ is the unit Euclidian ball. Notice that $\tilde{\mathcal{W}} = \{ w \in \mathbb{R}^d \; : \; w^\top w \leq  (1 + d)^2\}$ by the definition of Minkowski sum. At this point one can always define for any ${w} \in \hat{\mathcal{W}}$ a corresponding $\tilde{w} = (1+d) \tilde{w} \in \tilde{\mathcal{W}}$. By letting $\tilde{x} = (1+d) x$, it follows that $x^\top \Omega^{-1} x \leq 1$ if and only if $\tilde{x}^\top \tilde{\Omega}^{-1}\tilde{x} \leq 1$ with $\tilde{\Omega}^{-1} = (1+d)^2 \Omega^{-1}$.

Next, we move to the second case ($\hat{\mu} \neq \mu$). Assume the true measure $\mu$ is known to have zero mean, but its covariance $\Sigma$ is uncertain with only a nominal estimate $S \succ 0$ available. As a result, the nominal measure $x_{t+1}|x_t := \hat{\mu}_{x_{t+1}|x_t}$ with mean $Ax_t$ and covariance $S$ will inherit the misspecification. To account for this ambiguity, we revise \eqref{eq:safety:martingale} and require
\begin{equation}\label{eq:robust:cond}
\inf_{\nu \in \mathcal{A}(\hat{\mu}_{x_{t+1}|x_t})} \mathbb{E}[b(x_{t+1}) | \mathcal{F}_{t}]  \geq \tilde{\beta} b(x_t) + \delta,
\end{equation}
where $\mathcal{A}(\hat{\mu}_{x_{t+1}|x_t})$ is a judiciously chosen set of probability measures (i.e., a so-called \emph{ambiguity set}) encapsulating our trust in the nominal measure $\hat{\mu}$. Specifically, for $\rho > 0$, we define
\begin{equation}\label{eq:ambiguity:set}
\begin{aligned}
\mathcal{A}_\rho(\hat{\mu}_{x_{t+1}|x_t}) \coloneqq \{ & \nu \in \mathcal{P}_2(\mathbb{R}^d) \; :  \mathbb{E}_{\xi \sim \nu}[\xi] = 0,\\ &\mathbb{E}_{\xi \sim \nu}[\xi \xi^\top] = \Sigma \in  \mathcal{G}_\rho (S)\}
\end{aligned}
\end{equation}
where $\mathcal{P}_2(\mathbb{R}^d)$ is the space of probability measures in $\mathbb{R}^d$ with finite second moment,$$ \mathcal{G}_\rho (S) \coloneqq \left\{\Sigma \in  \mathbb{S}_+^n\; : \sqrt{\text{Tr}\left(\Sigma + S - 2\left(S^{\frac{1}{2}} \Sigma S^{\frac{1}{2}}\right)^{\frac{1}{2}}\right)} \leq \rho \right\}$$ 
is the co-called Gelbrich distance, and $\mathbb{S}_+^n$ is the set of symmetric positive definite matrices in $\mathbb{R}^{n \times n}$.

% \textit{Robust setting:}  Consider the set $\|\Sigma^{1/2} - S^{1/2}\|_F\leq c$ for $c > 0$ being a prescribed ambiguity budget describing our trust in the nominal probability distribution. To robustify the design in Subsection 3.2 against ambiguity in the parameter $\Sigma$, we  require the following implication to hold
% \begin{equation}\label{eq:robust:1}
% \|\Sigma^{1/2} - S^{1/2}\|_F \leq c \implies  \begin{bmatrix}
% \lambda I & \Sigma^{ 1/2} \\ \star & \Omega
% \end{bmatrix} \geq 0,
% \end{equation}
% in order to ensure \eqref{eq:safety:martingale} holds for all measures whose covariances belong to the prescribed set.
% By direct application of the S-Procedure (Lemma \ref{eq:S-procedure1}), \eqref{eq:robust:1} holds if and only if there exists a $\nu > 0$ such that
% \begin{equation}\label{eq:robust:sigma}
%  \begin{bmatrix}
% (\lambda - \nu) I & S \\ S^\top & \Omega - \nu I
% \end{bmatrix} \geq 0, \quad   \nu \geq \frac{c^2}{2}.
% \end{equation}
% Note that \eqref{eq:robust:sigma} consists of a set of convex constraints in $(\lambda, \nu, \Omega)$ that can be directly embedded into Problem \eqref{eq:convex-without-input-constraints-finite}.

\begin{proposition}\label{prop:gelbrich}
The distributionally robust constraint \eqref{eq:robust:cond} with ambiguity set as in \eqref{eq:ambiguity:set} is equivalent to the following feasibility problem: 
\begin{equation}\label{eq:gelbrich:ref}
\begin{aligned}
\emph{Find} \quad & \Omega, Y, \gamma, Z, Q, q_0\\
 \emph{s.t.} \quad &   q_0 + \gamma(\rho - \emph{Tr}(S))  + \emph{Tr}(Z) \leq 0\\  
 &  \gamma \geq 0, Z \in \mathbb{S}_+^n,  q_0 \geq (1-\tilde{\beta} + \delta), \\
 & \begin{bmatrix} Q & I\\ I & \Omega \end{bmatrix} \succeq 0,\: \begin{bmatrix} \gamma I - Q & \gamma S^{1/2} \\ \star & Z \end{bmatrix} \succeq 0\\
  &   q_0I \succeq \begin{bmatrix}
        \tilde{\beta}\Omega \:\: & \Omega A^\top\!\!+\!\!B^\top Y^\top \\
        \star&\Omega
    \end{bmatrix}.
\end{aligned}
\end{equation}
\end{proposition}
Eq. \eqref{eq:gelbrich:ref} represents a jointly convex program in all the decision variables, hence it can be embedded into \eqref{eq:convex-without-input-constraints-finite} in substitution of constraints (14d) and (14f).

\subsection{Safety filters}\label{sec:filter}
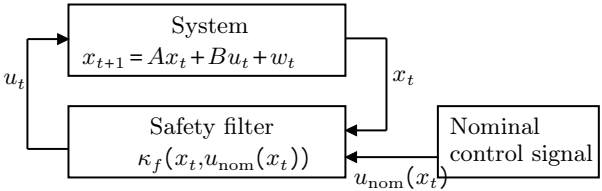
\begin{figure}[b!]
\centering

\tikzset{every picture/.style={line width=0.75pt}} %set default line width to 0.75pt        

\begin{tikzpicture}[x=0.75pt,y=0.75pt,yscale=-0.9,xscale=0.9]

%Shape: Rectangle [id:dp316490033936607] 
\draw   (109,40) -- (264.2,40) -- (264.2,80) -- (109,80) -- cycle ;
%Shape: Rectangle [id:dp9114000442540631] 
\draw   (109,97) -- (264.2,97) -- (264.2,137) -- (109,137) -- cycle ;
%Straight Lines [id:da2872102908842773] 
\draw    (264,59) -- (287.2,59) ;
%Straight Lines [id:da16671806591801874] 
\draw    (287.2,59) -- (287.2,110.6) ;
%Straight Lines [id:da8198610495730188] 
\draw    (287.2,110.6) -- (267.2,110.6) ;
\draw [shift={(264.2,110.6)}, rotate = 360] [fill={rgb, 255:red, 0; green, 0; blue, 0 }  ][line width=0.08]  [draw opacity=0] (6.25,-3) -- (0,0) -- (6.25,3) -- cycle    ;
%Straight Lines [id:da4936255814925381] 
\draw    (86,121) -- (109.2,121) ;
%Straight Lines [id:da48971005530540723] 
\draw    (86,59.4) -- (86,121) ;
%Straight Lines [id:da8496138121089776] 
\draw    (86,59.4) -- (106.2,59.4) ;
\draw [shift={(109.2,59.4)}, rotate = 180] [fill={rgb, 255:red, 0; green, 0; blue, 0 }  ][line width=0.08]  [draw opacity=0] (6.25,-3) -- (0,0) -- (6.25,3) -- cycle    ;
%Shape: Rectangle [id:dp8897907805977943] 
\draw   (316,97) -- (406.2,97) -- (406.2,137) -- (316,137) -- cycle ;
%Straight Lines [id:da4992261830006044] 
\draw    (316.2,125.6) -- (267.2,125.6) ;
\draw [shift={(264.2,125.6)}, rotate = 360] [fill={rgb, 255:red, 0; green, 0; blue, 0 }  ][line width=0.08]  [draw opacity=0] (6.25,-3) -- (0,0) -- (6.25,3) -- cycle    ;

% Text Node
\draw (165,45) node [anchor=north west][inner sep=0.75pt]  [font=\small] [align=left] {System};
% Text Node
\draw (115,62.4) node [anchor=north west][inner sep=0.75pt]  [font=\footnotesize]  {$x_{t+1}\! =\! Ax_{t} \! +\! Bu_{t} \! +\! w_{t}$};
% Text Node
\draw (153,102) node [anchor=north west][inner sep=0.75pt]  [font=\small] [align=left] {Safety filter};
% Text Node
\draw (146,118.4) node [anchor=north west][inner sep=0.75pt]  [font=\footnotesize]  {$\kappa _{f}( x_{t} ,\! u_{\text{nom}}( x_{t}))$};
% Text Node
\draw (71,77.4) node [anchor=north west][inner sep=0.75pt]  [font=\footnotesize]  {$u_{t}$};
% Text Node
\draw (289,75.4) node [anchor=north west][inner sep=0.75pt]  [font=\footnotesize]  {$x_{t}$};
% Text Node
\draw (321,102) node [anchor=north west][inner sep=0.75pt]  [font=\small] [align=left] {Nominal\\control signal};
% Text Node
\draw (267.2,129) node [anchor=north west][inner sep=0.75pt]  [font=\footnotesize]  {$u_{\text{nom}}( x_{t})$};

\end{tikzpicture}

\caption{Illustration of the safety filter concept.}
\label{fig:filter}
\end{figure}

Lastly, we point out that the construction presented in Section \ref{sec:main} can be employed to synthesize safety filters. To streamline the result, we again consider the infinite-horizon setting. A safety filter endows a well-performant input signal (often coming from stabilizing controllers
hand-designed by domain specialists, learning-based
controllers that maximize a particular reward signal, or
human input to the system with safety guarantees) with safety guarantees. More formally, a safety filter $\kappa_f : \mathbb{R}^n \times \mathbb{R}^m \rightarrow \mathbb{R}^m$ modifies
the nominal control input signal to produce an input signal
$u = \kappa_f ( x_t, u_{\text{nom}}(x_t))$ that ensures the system meets the safety requirement, while minimally modifying the
nominal input signal, i.e., minimizing
\begin{equation}\label{eq:cost:filter}
J( x_t, u_{\text{nom}}(x_t)) = \frac{1}{2}\|u -  u_{\text{nom}}(x_t)\|_2^2.
\end{equation}
A block diagram of a safety filter is shown in Figure \ref{fig:filter}.

Assume a valid CBF $b(x) = 1 - x^\top \Omega^{-1}x$ has been synthesized by solving \eqref{eq:convex-without-input-constraints}. Then, a safety filter can be realized by repeatedly solving a parametric optimization program in closed-loop for every state $x \in \mathbb{R}^n$ along the trajectory as
\begin{equation}\label{program:general}
\begin{array}{rll}
 \min _{u} & J( x_t, u_{\text{nom}}(x_t))  \\
\text { s.t. } & b(Ax_t + Bu + w) - b(x_t) \geq -\beta b(x_t), \forall w \in \mathcal{W},
\end{array}
\end{equation}
with $\beta \in (0,1)$ fixed. The objective of this section is to analyze the regularity properties of the minimizer of \eqref{program:general}, denoted as $u = \kappa_f ( x_t, u_{\text{nom}}(x_t))$, to ensure that the control action is sufficiently regular. We begin with the following instrumental result\footnote{The proof is a direct consequence of the convexity definition, hence it is omitted in the interest of space.}:
\begin{lemma}\label{lemma:convexity}
Let ${\beta}, x, w$ be given parameters. The constraint $b(Ax+Bu+w) \geq {(1-\beta)} b(x)$ is convex in $u$ whenever $b$ is concave.
\end{lemma}
% \begin{proof}
%  Consider $u_1, u_2$ satisfying the $b(x+u_1)\geq \beta b(x)$ and $b(x+u_2)\geq \beta b(x)$ for a fixed $x$, respectively. Take $u = \lambda u_1 + (1-\lambda) u_2$ for $\lambda \in [0,1]$. Then it holds:
% $$
% \begin{aligned}
% b(x + \lambda u_1 + (1-\lambda) u_2) & \overset{(a)}{\geq} \lambda b(x + u_1) + (1-\lambda)  b(x + u_2)\\
% & \overset{(b)}{\geq} \lambda(\beta b(x)) + (1-\lambda)(\beta b(x))\\
% & = \beta b(x),
% \end{aligned}
% $$
% where $(a)$ follows since $b$ is concave by assumption, and $(b)$ by the inequalities satisfied by $u_1$ and $u_2$.
% \end{proof}

By Lemma \ref{lemma:convexity}, \eqref{program:general} is a convex optimization program in $u$ for any fixed $x$ and $w$, as $b(x)$ synthesized with \eqref{eq:convex-without-input-constraints} is concave since $\Omega^{-1} \succ 0$. By exploiting the uncertainty description in Assumption \ref{eq:support} and the S-Lemma (Lemma \ref{eq:S-procedure1}), we can equivalently reformulate the semi-infinite problem in \eqref{program:general} as
\begin{equation}\label{program:reformulated}
\begin{array}{rll}
 \min _{u} & J( x_t, u_{\text{nom}}(x_t))  \\
\text { s.t. } & \lambda, \gamma \geq 0,\\
& u^\top Q u + 2u^\top a + \kappa + \gamma \leq 0,\\
& \begin{bmatrix} L + \lambda I & b u + d\\
\star & \gamma - \lambda\end{bmatrix} \succeq 0,
\end{array}
\end{equation}
where
$$
\begin{aligned}
Q & = B^\top \Omega^{-1}B,  &a & = B^\top \Omega^{-1}A x_t,\\
b & = D^\top \Omega^{-1}B, & L & = D^\top\Omega^{-1}D\\
d & = D^\top \Omega^{-1} A x_t,& \kappa & = x_t^\top(A^\top\Omega^{-1}A - \beta \Omega^{-1}) x_t + (\beta-1).
\end{aligned}
$$
Problem \ref{program:reformulated} is amenable to fast semidefinite programming solvers, such as \texttt{STROM} \cite{kang2024fast} which is shown to deliver real-time performance even for large-scale instances.
% We rely on the following assumption that is common in parametric nonlinear programming \cite{subotic2021quantitative} to argue about the properties of the designed safety filter.
% \begin{assumption}\label{ass:LICQ}
% We assume that Linear Independence Constraint Qualification (LICQ) is satisfied for all $t \geq 0$.
% \end{assumption}
% \begin{proposition}\label{prop:Lipschitz}
% Under Assumption \ref{ass:LICQ}, (1) problem \eqref{program:reformulated} is recursive feasible, and (2) the optimal solution of \eqref{program:reformulated} is locally Lipschitz continuous.
% \end{proposition}

\subsection{General safety specifications}

In Section \ref{sec:main}, the safe set $\mathcal{S}$ is restricted to be the intersection of half-planes. While this is satisfactory for many applications, it is not applicable for more complex safety specifications, such as obstacle avoidance. We extend our results to more general settings, where $\mathcal{S}$ and $\mathcal{I}$ are semi-algebraic sets defined in \eqref{eq:poly:set} and \eqref{eq:poly:set:init}, motivated by the following example.

\begin{example}[Robust motion-planning problem]
Consider a robust motion-planning problem for system \eqref{eq:poly:sys}. The task is to plan the motion of such system from an initial state $x_0 \in \mathcal{I} \subset\mathbb{R}^n$ to a target state $y_0 \in \mathcal{S}$ through an obstacle-filled environment $\bar{\mathcal{Y}} \subset \mathbb{R}^n$. The safe set $\mathcal{S}$ is modeled by the set difference of a nominal constraint set $\bar{\mathcal{Y}}$ and a collection of obstacles $\{\mathcal{O}_i\}_{i=1}^O$
$$
\mathcal{S} = \bar{\mathcal{Y}} \setminus \textstyle\left( \bigcup_{i,\ldots, O} \mathcal{O}_i \right).
$$
Assume that $\bar{\mathcal{Y}}$ and each $\mathcal{O}_i$ are semi-algebraic. Then, $\bigcup_{i,\ldots, O} \mathcal{O}_i$ is semi-algebraic as semi-algebraic sets are closed under finite union. The complement $\mathcal{R}^n \setminus \textstyle\left( \bigcup_{i,\ldots, O} \mathcal{O}_i \right)$ is also algebraic. Finally, the intersection $\bar{\mathcal{Y}} \bigcup (\mathcal{R}^n \setminus \textstyle\left( \bigcup_{i,\ldots, O} \mathcal{O}_i \right)) = \bar{\mathcal{Y}} \setminus \textstyle\left( \bigcup_{i,\ldots, O} \mathcal{O}_i \right)$ is semi-algebraic as the intersection of semi-algebraic sets remains algebraic.   $\hfill \square$
\end{example}

The problem of interest is again to design a control invariant set $\mathcal{B}$ for system \eqref{eq:poly:sys} under Assumption \ref{eq:support}, such that the set containment conditions $\mathcal{I}\subseteq \mathcal{B}$ and $\mathcal{B}\subseteq\mathcal{S}$ hold. Again for brevity we consider the setting from Section 3.1. We will rely on tools from SOS Programming, which we briefly recall next. 
\vspace{-0.2cm}
\begin{definition}[SOS Polynomial]\label{def:poly}
A polynomial $p(x)$ is said to be a sum-of-squares polynomial in $x \in \mathbb{R}^n$ if there exists $N$ polynomials $p_i(x), \:i = 1,\ldots,N$ such that
$
p(x) = \sum_{i=1}^N p_i(x)^2.
$
\end{definition}

% Clearly, if $p(x)$ is a sum-of-squares polynomial, then it is non-negative for all $x \in \mathbb{R}^n$. Further, computing the sum-of-squares decomposition in Def.~\ref{def:poly} amounts to a positive semidefinite feasibility program, which is convex.

\begin{lemma}\label{lemma:sos}
Consider a polynomial $p(x)$ of degree $2d$ in $x \in \mathbb{R}^n$ and let $z(x)$ be a vector of user-defined monomials of fixed degree less or equal than $d$. Then $p(x)$ admits a sum-of-squares decomposition if and only if
$$
p(x) = z(x)^\top M z(x), \: M \succeq 0.
$$
\end{lemma}
\vspace{-0.3cm}
For a given $z(x)$, the task of finding $M$ amounts to a semi-definite program, which can be solved efficiently
using interior point methods. Selecting the basis $z(x)$ depends on the structure of the polynomial $p(x)$ to be decomposed. 

We are now ready to outline the design program, which is shown in \eqref{eq:design-program}. 
\vspace{-0.1cm}
\begin{subequations}\label{eq:design-program}
    \begin{align}
        \mathrm{find}~&b(x),u(x)\mathrm{~s.t.~\eqref{eq:parametrization:all}},
        \sigma_b(x),\{\sigma_{s,j}(x)\}_{j=1}^q\in\mathrm{SOS}[x]\nonumber\\
        &\{\sigma_{r,j}(x)\}_{j=1}^p, \sigma_w(x)\in\mathrm{SOS}[x]\label{eq:decisionvariables}\\
        \mathrm{s.t.}~&b(Ax+Bu(x)+w)-\sigma_b(x)b(x)\nonumber\\
        &-\sigma_w(x)(1-w^\top w)\in\mathrm{SOS}[x,w]\label{eq:invariance-general}\\
        &b(x)-\sum_{j=1}^p\sigma_{r,j}(x)r_j(x)\in\mathrm{SOS}[x],\label{eq:set2}\\
        &-b(x)+\sigma_{s,j}(x)s_j(x)\in\mathrm{SOS},j\in\{1,\ldots,q\}.\label{eq:set3}
    \end{align}
\end{subequations}
\vspace{-0.3cm}

Note that the constraint \eqref{eq:invariance-general} is nonconvex due to the term $b(Ax+Bu(x)+w)$. This also prevents the applicability of bilinear SDP solvers such as PENLAB \cite{fiala2013penlab}, as well as iterative methods such $\mathcal{KL}$-iteration for the problem. In Section \ref{sec:main}, this constraint is effectively convexified by changing variables and lifting the dimension of the matrix inequality. The set containment constraints can be convexified simultaneously when $\mathcal{S}$ and $\mathcal{I}$ take special forms. However, these techniques can not be directly applied to convexify the program \eqref{eq:design-program}, as the set containment constraints \eqref{eq:set2} and \eqref{eq:set3} include higher order polynomials, thus remaining nonconvex. We propose Algorithm \ref{al:nonlinear} to solve \eqref{eq:design-program} in an iterative\footnote{We use $\Omega_k$ and $Y_k$ to denote the solution of $\Omega$ and $Y$ at the $k$-th iteration.}manner with guaranteed convergence to a feasible CBF $b(x)$ and controller $u(x)$.

\begin{algorithm}
 \caption{Iterative algorithm for \eqref{eq:design-program}}
  \hspace*{\algorithmicindent} \textbf{Initialization} Tolerance $\texttt{tol}$, $\Omega_0$ as the solution of \eqref{eq:initialization}:
  \begin{subequations}\label{eq:initialization}
\begin{align}
    \min~&\mathrm{Tr}(\Omega)\\
    \mathrm{s.t.}~&  Y \in\mathbb{R}^{m\times n},  0\prec\Omega=\Omega^\top\in\mathbb{R}^{n\times n},  0\prec R\in\mathbb{R}^{n\times n},\\ 
   &  \{\sigma_{r,j}(x)\}_{j=1}^p\in\Sigma[x],\\
    &\begin{bmatrix}
        -\tilde{\beta} \Omega + \lambda \Omega & 0 & \Omega A^\top + Y^\top B^\top \\
        0 & -\lambda I & D^\top\\
        A\Omega + BY & D & -\Omega
    \end{bmatrix} \preceq 0,\label{eq:invariance-initialization}\\
     &\begin{bmatrix}
     R&I\\I&\Omega
     \end{bmatrix}\succeq 0,\label{eq:containment1}\\
        &1-x^\top R x-\sum_{j=1}^p\sigma_{r,j}(x)r_j(x)\in\mathrm{SOS}[x].\label{eq:hahaha}
\end{align}
\end{subequations}\\
 \hspace*{\algorithmicindent} Let $b_0(x)=1-x^\top \Omega_0^{-1}x.$\\
 \vspace{-3ex}
 \begin{algorithmic}[1]\label{al:nonlinear}
 \WHILE{$\|\Omega_{k}-\Omega_{k-1}\|_F>\texttt{tol}$}
 \STATE \label{step:firstcomputation} \textbf{solve} $Y_k$ as the solution of the following optimization problem
 \begin{equation}\label{eq:prog1}
   \begin{split}
     \mathop{\mathrm{find}}_{Y\in\mathbb{R}^{m\times n}}~&Y\\
     \mathrm{s.t.~}~&\begin{bmatrix}
        -\tilde{\beta} \Omega_{k-1} \!+\! \lambda \Omega_{k-1} & 0 & \Omega_{k-1} A^\top \!+\! Y^\top B^\top \\
        0 & -\lambda I & D^\top\\
        A\Omega_{k-1} \!+ \!BY & D & -\Omega_{k-1}
    \end{bmatrix} \!\preceq \!0
   \end{split}
 \end{equation}
 Let $u_k(x)=Y_k\Omega_{k-1}^{-1}x$.
 \STATE \textbf{solve} $\Omega_k'$ as the optimal solution of $\Omega$ for the following optimization problem in decision variables $\Omega$ and $\{\sigma_{\{\cdot\}}(x)\}~\mathrm{as~in~\eqref{eq:decisionvariables}}$.
 \begin{equation}\label{eq:prog2}
     \begin{split}
         \min_{b(x)=1-x^\top \Omega x,u(x)=u_k(x)}~&\Omega\\
         \mathrm{subject~to~}&\eqref{eq:invariance-general}-\eqref{eq:set3}\\
         &
     \end{split}
 \end{equation}
 Let $\Omega_k=\Omega_k'^{-1}$.
 \ENDWHILE
 \end{algorithmic}
\end{algorithm}

Algorithm \ref{al:nonlinear} uses the optimal solution of \eqref{eq:initialization} to initialize the CBF $b_0(x)=1-x^\top \Omega_0^{-1}x$. Constraint \eqref{eq:invariance-initialization} corresponds to the control invariance condition. Constraint \eqref{eq:containment1} ensured $\mathcal{R}:=\{x\in\mathbb{R}^n:1-x^\top Rx\ge 0\}\subseteq \mathcal{B}:=\{x\in\mathbb{R}^n:1-x^\top \Omega^{-1}x\ge 0\}$. Constraint \eqref{eq:hahaha} is sufficient for $\mathcal{I}\subseteq \mathcal{R}$. As a result, the CBF $b_0(x)=1-x^\top \Omega_0^{-1}x$ defines a control invariant set $\mathcal{B}$ that satisfies $\mathcal{I}\subseteq\mathcal{B}$. The objective function is to minimize the volume of $\mathcal{B}$. If $\mathcal{B}\subseteq \mathcal{S}$ holds, then $b_0(x)$ is a feasible initialization that warm-starts Algorithm \ref{al:nonlinear}. However, this is not ensured as the nonconvex set containment condition $\mathcal{B}\subseteq \mathcal{S}$ is not incorporated into \eqref{eq:initialization}.

After obtaining the initial CBF $b_0(x)$, we iteratively solve \eqref{eq:prog1} for control $u_k(x)$, and apply the control into \eqref{eq:prog2} for CBF $b_k(x)$. It should be noted that the CBF $b(x)$ in \eqref{eq:prog2} has been parameterized by $b(x)=1-x^\top \Omega x$ but not $1-x^\top \Omega^{-1}x$ as in \eqref{eq:parametrization:all}. The new parameterization ensures convexify of constraints \eqref{eq:set2}-\eqref{eq:set3}. Moreover, the control invariance constraint \eqref{eq:invariance-general} is also convex without changing variables and lifting dimension of the matrix inequality. This is because the control input $u(x)$ has been replaced by $Y_k\Omega_{k-1}^{-1}x$, where $Y_k$ and $\Omega_{k-1}$ are known. Then, $b(Ax+Bu_k(x))$ is inherently linear in $\Omega$, demonstrating convexity of constraint \eqref{eq:invariance-general}. Similarly, constraints \eqref{eq:set2} - \eqref{eq:set3} are all linear in $\Omega$. After solving \eqref{eq:prog2} for $\Omega_{k}'$, we calculate $\Omega_k=\Omega_{k-1}'^{-1}$, so that $b(x)=1-x^\top \Omega_k^{-1}x$.

\section{Numerical examples}\label{sec:5}
\subsubsection*{Example 1 - Pendulum} 
\begin{figure*}
\centering
\begin{minipage}[t]{0.3\textwidth}
    \centering
    \includegraphics[width=1.16\linewidth]{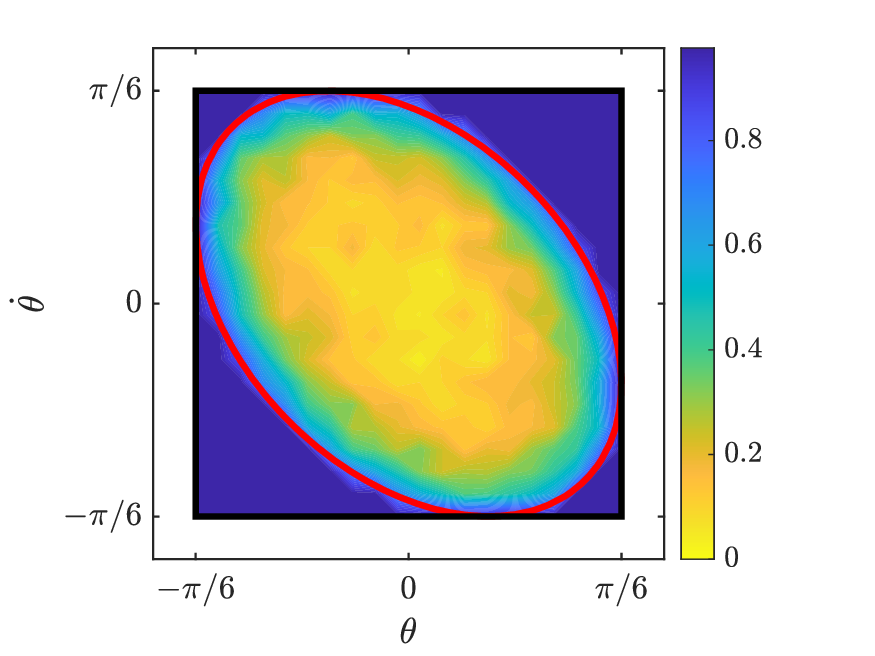}
    \caption{Empirical bound.}\label{fig:pendulum:sensitivity:emp}
\end{minipage}
\begin{minipage}[t]{0.3\textwidth}
    \centering
    \includegraphics[width=1.09\linewidth]{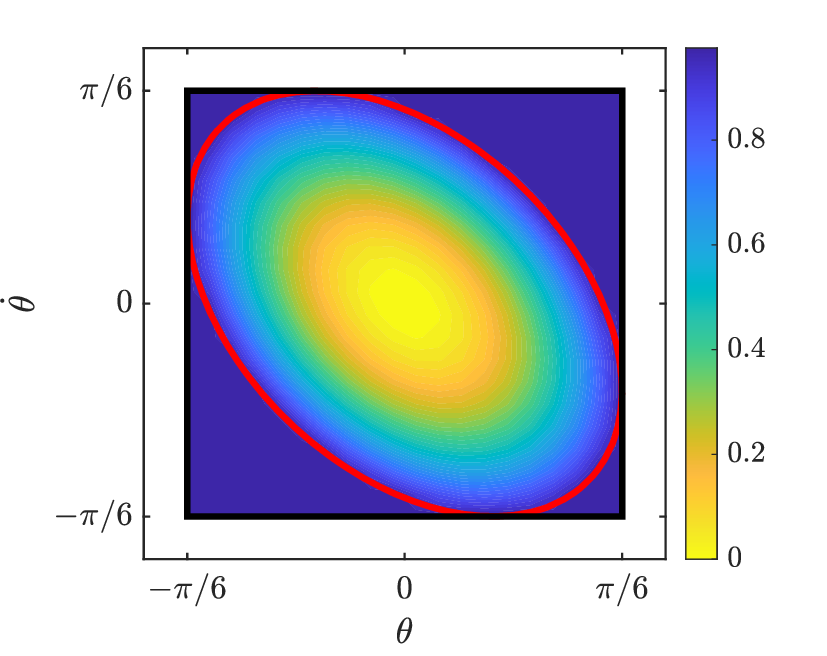}
    \caption{Theoretical bound.}\label{fig:pendulum:sensitivity}
\end{minipage}
\begin{minipage}[t]{0.39\textwidth}
    \centering
    \includegraphics[width=1\linewidth]{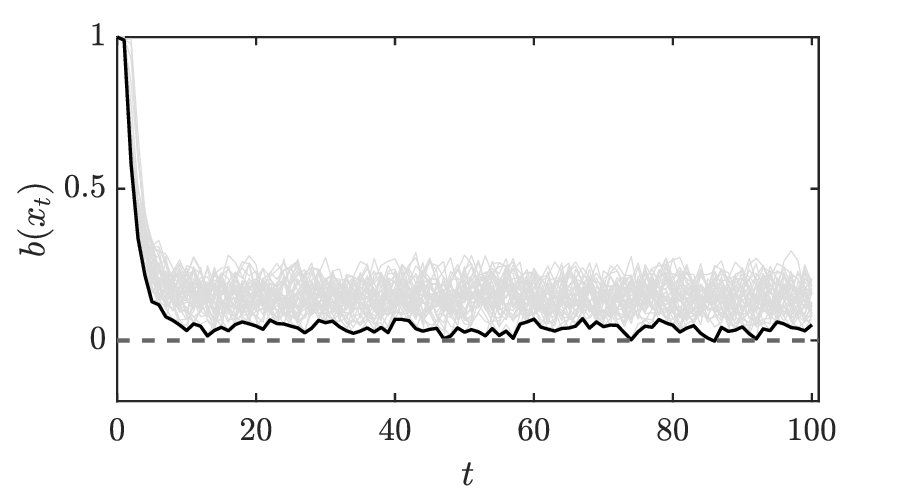}
    \caption{100 one-second-long trajectories.}\label{fig:pendulum:trajectories}
\end{minipage}%
\end{figure*}

We consider an inverted pendulum to be stabilized about its upright
equilibrium point with the following discrete-time dynamics linearized around the origin:
% $$
% \begin{bmatrix}  \theta_{t+1} \\ \dot{\theta}_{t+1} \end{bmatrix} = \begin{bmatrix}  \theta_t + \Delta t \dot{\theta}_t\\ \dot{\theta}_t + \Delta t \text{sin}(\theta_t)\end{bmatrix} + \begin{bmatrix} 0\\  \Delta t u\end{bmatrix} + w_t,
% $$
% We linearized the above dynamics around the origin obtaining 
$$
\begin{bmatrix}  \theta_{t+1} \\ \dot{\theta}_{t+1} \end{bmatrix} = \begin{bmatrix}  1 & \Delta t\\ \Delta t & 1\end{bmatrix}\begin{bmatrix}  \theta_{t} \\ \dot{\theta}_{t} \end{bmatrix} + \begin{bmatrix} 0\\  \Delta t \end{bmatrix}u_t + w_t,
$$
with step size $\Delta t = 0.01$ s, and i.i.d. disturbances $w_t \sim \mathcal{N}(0_2, \text{diag}(0.0075^2, 0.05^2))$. This linearization is a good approximation of the original nonlinear system for $|\theta| \leq \frac{\pi}{6}, |\dot{\theta}| \leq \frac{\pi}{6}$. As a result, we consider as a safe set to employ our linearization, the rectangle described by
$$
\mathcal{S} = \left\{ \begin{bmatrix}  \theta_{t+1} \\ \dot{\theta}_{t+1} \end{bmatrix} \in \mathbb{R}^2 \; : \; \begin{bmatrix}  \theta_{t+1} \\ \dot{\theta}_{t+1} \end{bmatrix} \in \left[-\frac{\pi}{6},\frac{\pi}{6}\right] \times \left[-\frac{\pi}{6},\frac{\pi}{6}\right] \right\}.
$$
We consider a simulation time of $1$ s, resulting in $T = 100$ time steps. We construct a CBF $b(x)$ and a feedback controller $u(x)$ by solving Problem \eqref{eq:convex-without-input-constraints-finite} with fixed $\beta = 0.8$ and $\delta$ = 0. Note that this construction avoids the \textit{manual} derivation of a CBF as done in \cite{taylor2022safety} based on the corresponding continuous-time Lyapunov equation.  Figure \eqref{fig:pendulum:trajectories} shows the result of 500 Monte-Carlo simulations starting at $x_0 = [0,0]^\top$, with an empirical safety level of $91\%$. Finally, we perform a sensitivity analysis to understand the effect of the initial condition on the bound in Theorem \eqref{th:convex-design}. To this end, we grid the state space with equally spaced points: for each gridded $x_0$ we run 100 Monte-Carlo simulations and we compute the theoretical bound resulting from Proposition \ref{prop:Ville:general}. Figure \eqref{fig:pendulum:sensitivity} compares the empirical result (left) with the theoretical result (right). As expected from Proposition \ref{prop:Ville:general}, the safer the initial condition, the lower the exit probability (hence, the higher the safety probability). Additionally, the two subplots highlight the conservativism of the theoretical bound (which in fact only provides a lower bound to the safety probability). 

% \begin{figure}
%     \centering
%     \includegraphics[width=0.8\linewidth]{graphics/stochastic_traj.eps}
%     \caption{00 one-second-long trajectories (in grey) of the pendulum starting at the origin. In red we highlight an example of a violating trajectory, while in green we highlight a non-violating one. }
%     \label{fig:pendulum:trajectories}
% \end{figure}

% \begin{figure}
%     \centering
%     \includegraphics[width=\linewidth]{graphics/stochastic_plot.eps}
%     \caption{Left: Montecarlo estimate of the exit probability for different initial conditions. Right: theoretical bound from Proposition \ref{prop:Ville}. Safety probability is computed wrt $\mathcal{B}$, hence for $\mathcal{S}$ is in general higher.}
%     \label{fig:pendulum:sensitivity}
% \end{figure}

\subsubsection*{Example 2 - Double integrator} We consider a double-integrator with dynamics
$$
\begin{bmatrix}  x^1_{t+1} \\ x^2_{t+1} \end{bmatrix} = \begin{bmatrix}  0.1 & 0.65 \\ 0 & 1.02\end{bmatrix} \begin{bmatrix}  x^1_{t} \\ x^2_{t} \end{bmatrix} + \begin{bmatrix} 0.5 \\0.5\end{bmatrix}u_t +  \begin{bmatrix} 0.01 & 0 \\0 & 0.01\end{bmatrix} w_t,
$$
with $w \in \mathbb{W}$ as defined in Assumption \eqref{eq:support}. The system is equipped with a nominal controller $u_{\text{nom}}(x_t) = [0, 50] x_t$ that steers the system towards the upper-right direction that maximizes a given performance score. However, due to safety reasons, the system states are required to lie within the safe set
$$
\mathcal{S} = \left\{ \begin{bmatrix}  x^1_{t} \\ x^2_{t} \end{bmatrix} \in \mathbb{R}^2 \; : \; \begin{bmatrix}  x^1_{t} \\ x^2_{t} \end{bmatrix} \in \left[ -2,2\right] \times \left[ -2,2\right] \right\}.
$$
As a result, we design a safety filter to enhance the safety of the nominal controller that solves in real-time Problem \ref{program:reformulated}, following the procedure described in Subsection \ref{sec:filter}. We consider 50 Montecarlo simulations of length $T = 100$ starting from the origin. We use $\beta = 0.4$ and $\lambda = 0.05$. Figure \ref{fig:pendulum:trajectories} shows the evolution of $h(x_t)$ over time for each trajectory, together with the minimum value registered across the entire simulation campaign. As $h(x_t)\geq 0$, we conclude that the safety filter is capable of maintaining all the trajectories within the invariant set $\mathcal{B}$ and consequently within the safe set $\mathcal{S}$. 

% \begin{figure}
%     \centering
%     \includegraphics[width=0.9\linewidth]{graphics/example2_robust_small.eps}
%     \caption{Simulation results for the double integrator over 50 trials with safety filter in place. The grey curve display the value of the CBF across time for each trajectory, while the black solid line represents the minimum value of $b(x_t)$ register across time.}
%     \label{fig:integrator}
% \end{figure}

\section{Conclusions}\label{sec:6}
We considered discrete-time uncertain linear systems subject to additive disturbance and devised convex design rules for the joint synthesis of a control barrier function and a linear feedback law that guarantees safety. In particular, depending on the characteristics of the support set of the disturbance distribution, we considered two notions of safety: infinite-horizon and (joint-in-time) finite-horizon. We further extended our design procedure to account for several cases of practical relevance. Future work may include the extension of the co-design procedure to nonlinear systems via SOS relaxations and the design of predictive safety filters inspired by \cite{wabersich2022predictive} to incorporate anticipation of future states.

%%%%%%%%%%%%%%%%%%%%%%%%%%%%%%%%%%%%%%%%%%%%%%%%%%%%%%%%%%%%%%%%%%%%%%%%%%%%%%%%%%%%%%%%%%%%%%%%%%%%%%%%%%%%%%%%%%%%%%

%%%%%%%%%%%%%%%%%%%%%%%%%%%%%%%%%%%%%%%%%%%%%%%%

\bibliographystyle{plain}        % Include this if you use bibtex 
\bibliography{references}           % and a bib file to produce the 

\begin{thebibliography}{10}

\bibitem{ames2014control}
Aaron~D Ames, Jessy~W Grizzle, and Paulo Tabuada.
\newblock Control barrier function based quadratic programs with application to adaptive cruise control.
\newblock In {\em 53rd IEEE Conference on Decision and Control}, pages 6271--6278. IEEE, 2014.

\bibitem{bertsekas1972infinite}
Dimitri Bertsekas.
\newblock Infinite time reachability of state-space regions by using feedback control.
\newblock {\em IEEE Transactions on Automatic Control}, 17(5):604--613, 1972.

\bibitem{bojarski2016end}
Mariusz Bojarski.
\newblock End to end learning for self-driving cars.
\newblock {\em arXiv preprint arXiv:1604.07316}, 2016.

\bibitem{clark2022semi}
Andrew Clark.
\newblock A semi-algebraic framework for verification and synthesis of control barrier functions.
\newblock {\em arXiv preprint arXiv:2209.00081}, 2022.

\bibitem{cosner2023robust}
Ryan~K Cosner, Preston Culbertson, Andrew~J Taylor, and Aaron~D Ames.
\newblock Robust safety under stochastic uncertainty with discrete-time control barrier functions.
\newblock {\em arXiv preprint arXiv:2302.07469}, 2023.

\bibitem{dawson2023safe}
Charles Dawson, Sicun Gao, and Chuchu Fan.
\newblock Safe control with learned certificates: A survey of neural lyapunov, barrier, and contraction methods for robotics and control.
\newblock {\em IEEE Transactions on Robotics}, 39(3):1749--1767, 2023.

\bibitem{fiala2013penlab}
Jan Fiala, Michal Ko{\v{c}}vara, and Michael Stingl.
\newblock Penlab: A matlab solver for nonlinear semidefinite optimization.
\newblock {\em arXiv preprint arXiv:1311.5240}, 2013.

\bibitem{jagtap2020formal}
Pushpak Jagtap, Sadegh Soudjani, and Majid Zamani.
\newblock Formal synthesis of stochastic systems via control barrier certificates.
\newblock {\em IEEE Transactions on Automatic Control}, 66(7):3097--3110, 2020.

\bibitem{kang2024fast}
Shucheng Kang, Xiaoyang Xu, Jay Sarva, Ling Liang, and Heng Yang.
\newblock Fast and certifiable trajectory optimization.
\newblock {\em arXiv preprint arXiv:2406.05846}, 2024.

\bibitem{korda2014convex}
Milan Korda, Didier Henrion, and Colin~N Jones.
\newblock Convex computation of the maximum controlled invariant set for polynomial control systems.
\newblock {\em SIAM Journal on Control and Optimization}, 52(5):2944--2969, 2014.

\bibitem{kushner1967stochastic}
Harold~Joseph Kushner.
\newblock Stochastic stability and control.
\newblock 1967.

\bibitem{long2024distributionally}
Kehan Long, Jorge Cortes, and Nikolay Atanasov.
\newblock Distributionally robust policy and lyapunov-certificate learning.
\newblock {\em arXiv preprint arXiv:2404.03017}, 2024.

\bibitem{margellos2011hamilton}
Kostas Margellos and John Lygeros.
\newblock Hamilton--jacobi formulation for reach--avoid differential games.
\newblock {\em IEEE Transactions on automatic control}, 56(8):1849--1861, 2011.

\bibitem{masarotto2019procrustes}
Valentina Masarotto, Victor~M Panaretos, and Yoav Zemel.
\newblock Procrustes metrics on covariance operators and optimal transportation of gaussian processes.
\newblock {\em Sankhya A}, 81:172--213, 2019.

\bibitem{massiani2022safe}
Pierre-Fran{\c{c}}ois Massiani, Steve Heim, Friedrich Solowjow, and Sebastian Trimpe.
\newblock Safe value functions.
\newblock {\em IEEE Transactions on Automatic Control}, 68(5):2743--2757, 2022.

\bibitem{mitchell2005time}
Ian~M Mitchell, Alexandre~M Bayen, and Claire~J Tomlin.
\newblock A time-dependent hamilton-jacobi formulation of reachable sets for continuous dynamic games.
\newblock {\em IEEE Transactions on automatic control}, 50(7):947--957, 2005.

\bibitem{nguyen2021mean}
Viet~Anh Nguyen, Soroosh Shafiee, Damir Filipovi{\'c}, and Daniel Kuhn.
\newblock Mean-covariance robust risk measurement.
\newblock {\em arXiv preprint arXiv:2112.09959}, 2021.

\bibitem{parrilo2003semidefinite}
Pablo~A Parrilo.
\newblock Semidefinite programming relaxations for semialgebraic problems.
\newblock {\em Mathematical programming}, 96:293--320, 2003.

\bibitem{prajna2006barrier}
Stephen Prajna.
\newblock Barrier certificates for nonlinear model validation.
\newblock {\em Automatica}, 42(1):117--126, 2006.

\bibitem{prajna2007framework}
Stephen Prajna, Ali Jadbabaie, and George~J Pappas.
\newblock A framework for worst-case and stochastic safety verification using barrier certificates.
\newblock {\em IEEE Transactions on Automatic Control}, 52(8):1415--1428, 2007.

\bibitem{prandini2008application}
Maria Prandini and Jianghai Hu.
\newblock Application of reachability analysis for stochastic hybrid systems to aircraft conflict prediction.
\newblock In {\em 2008 47th IEEE conference on decision and control}, pages 4036--4041. IEEE, 2008.

\bibitem{reher2021dynamic}
Jenna Reher and Aaron~D Ames.
\newblock Dynamic walking: Toward agile and efficient bipedal robots.
\newblock {\em Annual Review of Control, Robotics, and Autonomous Systems}, 4:535--572, 2021.

\bibitem{robey2020learning}
Alexander Robey, Haimin Hu, Lars Lindemann, Hanwen Zhang, Dimos~V Dimarogonas, Stephen Tu, and Nikolai Matni.
\newblock Learning control barrier functions from expert demonstrations.
\newblock In {\em 2020 59th IEEE Conference on Decision and Control (CDC)}, pages 3717--3724. IEEE, 2020.

\bibitem{rockafellar2009variational}
R~Tyrrell Rockafellar and Roger J-B Wets.
\newblock {\em Variational analysis}, volume 317.
\newblock Springer Science \& Business Media, 2009.

\bibitem{schneeberger2024advanced}
Michael Schneeberger, Silvia Mastellone, and Florian D{\"o}rfler.
\newblock Advanced safety filter based on sos control barrier and lyapunov functions.
\newblock {\em arXiv preprint arXiv:2401.06901}, 2024.

\bibitem{singletary2022safe}
Andrew Singletary, Mohamadreza Ahmadi, and Aaron~D Ames.
\newblock Safe control for nonlinear systems with stochastic uncertainty via risk control barrier functions.
\newblock {\em IEEE Control Systems Letters}, 7:349--354, 2022.

\bibitem{srinivasan2020synthesis}
Mohit Srinivasan, Amogh Dabholkar, Samuel Coogan, and Patricio~A Vela.
\newblock Synthesis of control barrier functions using a supervised machine learning approach.
\newblock In {\em 2020 IEEE/RSJ International Conference on Intelligent Robots and Systems (IROS)}, pages 7139--7145. IEEE, 2020.

\bibitem{tan2004searching}
Weehong Tan and Andrew Packard.
\newblock Searching for control lyapunov functions using sums of squares programming.
\newblock {\em sibi}, 1(1), 2004.

\bibitem{taylor2022safety}
Andrew~J Taylor, Victor~D Dorobantu, Ryan~K Cosner, Yisong Yue, and Aaron~D Ames.
\newblock Safety of sampled-data systems with control barrier functions via approximate discrete time models.
\newblock In {\em 2022 IEEE 61st Conference on Decision and Control (CDC)}, pages 7127--7134. IEEE, 2022.

\bibitem{ville1939etude}
Jean Ville.
\newblock {\em Etude critique de la notion de collectif}.
\newblock Gauthier-Villars Paris, 1939.

\bibitem{wabersich2022predictive}
Kim~P Wabersich and Melanie~N Zeilinger.
\newblock Predictive control barrier functions: Enhanced safety mechanisms for learning-based control.
\newblock {\em IEEE Transactions on Automatic Control}, 68(5):2638--2651, 2022.

\bibitem{wang2023assessing}
Han Wang, Kostas Margellos, and Antonis Papachristodoulou.
\newblock Assessing safety for control systems using sum-of-squares programming, 2023.

\bibitem{wang2023safety}
Han Wang, Kostas Margellos, and Antonis Papachristodoulou.
\newblock Safety verification and controller synthesis for systems with input constraints.
\newblock {\em IFAC-PapersOnLine}, 56(2):1698--1703, 2023.

\bibitem{wang2024convex}
Han Wang, Kostas Margellos, Antonis Papachristodoulou, and Claudio De~Persis.
\newblock Convex co-design of control barrier function and safe feedback controller under input constraints.
\newblock {\em arXiv preprint arXiv:2403.11763}, 2024.

\bibitem{wang2024simultaneous}
Xinyu Wang, Luzia Knoedler, Frederik~Baymler Mathiesen, and Javier Alonso-Mora.
\newblock Simultaneous synthesis and verification of neural control barrier functions through branch-and-bound verification-in-the-loop training.
\newblock In {\em 2024 European Control Conference (ECC)}, pages 571--578. IEEE, 2024.

\bibitem{yang2023synthesizing}
Yujie Yang, Yuhang Zhang, Wenjun Zou, Jianyu Chen, Yuming Yin, and Shengbo~Eben Li.
\newblock Synthesizing control barrier functions with feasible region iteration for safe reinforcement learning.
\newblock {\em IEEE Transactions on Automatic Control}, 2023.

\bibitem{zhao2023convex}
Pan Zhao, Reza Ghabcheloo, Yikun Cheng, Hossein Abdi, and Naira Hovakimyan.
\newblock Convex synthesis of control barrier functions under input constraints.
\newblock {\em IEEE Control Systems Letters}, 2023.

\end{thebibliography}
                                 % bibliography (preferred). The
                                 % correct style is generated by
                                 % Elsevier at the time of printing.

%%%%%%%%%%%%%%%%%%%%%%%%%%%%%%%%%%%%%%%%%%%%%%%%

\section*{Appendix A: Proofs of Section 3}

\subsection{Proof of Theorem 3.3}
 We first prove that \eqref{eq:convex-e} is sufficient for $\mathcal{B} = \{ x\in \mathbb{R}^n \; : \; 1 - x^\top \Omega^{-1}x \geq 0\}$ to be an invariant set for $x_{t+1} = Ax_t + Bu(x_t) + Dw_t$. To this end, recall the robust forward invariant condition in \eqref{eq:robust:c2} and define the function $V(x_t) = x_t^\top \Omega^{-1}x_t$. Given the parametrization in \eqref{eq:parametrization:all}, define
 $$
 \begin{aligned}
 \Delta V & := V(x_{t+1}) - V(x_t) \\
 & = x_{t+1}^\top \Omega^{-1}x_{t+1} - \tilde{\beta}(x_t^\top \Omega^{-1}x_t)\\
 & = ((A + BY\Omega^{-1})x + Dw)^\top\Omega^{-1}((A + BY\Omega^{-1})x + Dw) \\
  & \quad - \tilde{\beta}(x_t^\top \Omega^{-1}x_t).
 \end{aligned}
 $$ 
Then, the right-hand side of \eqref{eq:robust:c2} is implied by the condition $\Delta V \leq 0$ since $\beta \in (0,1)$. We next show that this condition holds true for all $w \in \mathcal{W}$. Let $y = \begin{bmatrix} x^\top & w^\top &\left(\tilde{A}x + D w\right)^\top\end{bmatrix}^\top$ where $\tilde{A} = (A+BY\Omega^{-1})$ is the closed-loop matrix. Further, define
 $$
 M = \begin{bmatrix} -\tilde{\beta} \Omega^{-1} + \lambda \Omega^{-1} & 0 & \tilde{A}^\top \Omega^{-1}\\
 0 & -\lambda I & D^\top \Omega^{-1}\\
 \Omega^{-1}\tilde{A} & \Omega^{-1} D & -\Omega^{-1}\end{bmatrix}.
 $$
 Let $M  \preceq 0$ and notice that by post- and pre-multiplying $M$ by $\text{diag}(\Omega, I, \Omega)$ we exactly retrieve the condition in \eqref{eq:convex-e}. Next, since $M \preceq 0$ it holds $y^\top M y \leq 0$. In turn\footnote{With a slight abuse of notation, we substitute $x_t$ with $x$ and neglect temporal series relationships when clear from the context.}
 $$
 \begin{aligned}
 & y^\top M y  =   -x^\top \tilde{\beta} \Omega^{-1} x + x^\top \lambda \Omega x + x^\top \tilde{A}^\top \Omega \tilde{A} x +\\
 &  2 x^\top \tilde{A}^\top \Omega^{-1}Dw - \lambda w^\top w + w^\top D^\top \Omega^{-1} D w)\\ 
=  & x^\top (A+BY\Omega^{-1})^\top \Omega^{-1} (A+BY\Omega^{-1}) x + \\ 
 & 2 w^\top D^\top \Omega^{-1} (A+BY\Omega^{-1}) x + w^\top D^\top \Omega^{-1} D w - \tilde{\beta} x^\top  \Omega^{-1} x \\
 \leq & \lambda (w^\top w - x^\top \Omega^{-1} x).
 \end{aligned}
 $$
 By definition of $\Delta V$ and considering the control policy in \eqref{eq:parametrization:all}, this is equivalent to 
 $$
 \Delta V\leq \lambda (w^\top w - x^\top \Omega^{-1} x).
 $$
Recall now that $\lambda \in (0,1)$ and $w^\top w \leq 1$ by Assumption \ref{eq:support}. We distinguish two cases. For $x^\top \Omega^{-1} x > 1$, corresponding to $x \notin \mathcal{B}$, it holds 
$$ \Delta V \leq \lambda (w^\top w - x^\top \Omega^{-1} x) \leq 0.$$ 
Next, we consider the case $x^\top \Omega^{-1} x \leq 1$ corresponding to $x \in \mathcal{B}$. Notice that by definition of $V(x)$ we have
$$
\begin{aligned}
V(x_{t+1})
\leq & \lambda w^\top w + (1-\lambda) x^\top \Omega^{-1} x\\
\leq & \lambda + (1-\lambda) x^\top \Omega^{-1} x\\
\leq & \lambda + (1-\lambda) = 1,
\end{aligned}
$$
where the second inequality follows by Assumption \ref{eq:support} and the third one from the assumption $x^\top \Omega^{-1} x \leq 1$. In turn, $V(x_{t+1}) \leq 1$ implies $b(x_{t+1}) \geq 0$. 

We then consider $\mathcal{I}\subseteq \mathcal{B}$. From \eqref{eq:convex-f}, we have $R-\Omega^{-1}\succ 0$. As a result, we have
    \begin{equation*}
        \begin{split}
            \{x\in\mathbb{R}^n:1-x^\top Rx\ge 0\}\subseteq \{x\in\mathbb{R}^n:1-x^\top \Omega^{-1}x\ge 0\}.
        \end{split}
    \end{equation*}
    Therefore, we prove that $\mathcal{I}\subseteq\mathcal{B}$. We finally prove that \eqref{eq:convex-h} imply $\mathcal{B}\subseteq \mathcal{S}$. Using S-lemma \ref{eq:S-procedure1} for affine functions $a_j^\top x+1,j=1,\ldots,q$, and convex quadratic function $1-x^\top \Omega^{-1}x$, we have that $\mathcal{B}\subseteq\mathcal{S}$ if and only if for every $j=1,\ldots,q$, there exists $\lambda_i\ge \frac{1}{2}$ such that
\begin{equation*}
    2\lambda_j(a_j^\top x +1)-(-x^\top \Omega^{-1}x+1)\ge 0,\forall x\in\mathbb{R}^n,
\end{equation*}
which is equivalent to
\begin{equation}\label{eq:16}
    \forall j=1,\ldots,o,\exists \lambda_j\ge \frac{1}{2},\mathrm{s.t.}\begin{bmatrix}
        \Omega^{-1}&\lambda_ja_j\\
        \lambda_ja_j^\top & 2\lambda_j-1
    \end{bmatrix}\succeq 0.
\end{equation}
By Schur complement, \eqref{eq:16} holds if and only if $\Omega\succ 0$ and there exists $\lambda_j\ge \frac{1}{2}$ such that $2\lambda_j-1-\lambda^2a_j^\top \Omega a_j\ge 0$. The discriminant of the quadratic polynomial on the left-hand side of the inequality is $4-4a_j^\top \Omega a_j$. Hence, there exists $\lambda_j$ such that $2\lambda_i-1-\lambda_j^2a_j^\top \Omega a_i\ge 0$ if and only if $1-a_i^\top \Omega a_i\ge 0$. Moreover, if $1-a_i^\top \Omega a_i=0$, then $\lambda_i=1\ge \frac{1}{2},$ and if $1-a_i^\top \Omega a_i>0$, then any $\lambda_i\in[1-\sqrt{1-a_i^\top \Omega a_i},1+\sqrt{1-a_i^\top \Omega a_i}]$ satisfies $2\lambda_i-1-\lambda_i^2a_i^\top \Omega a_i\ge 0$. As $1+\sqrt{1-a_i^\top \Omega a_i}>\frac{1}{2},$ we have shown that there exists $\lambda_i\ge\frac{1}{2}$ such that $2\lambda_i-1-\lambda_i^2a_i^\top \Omega a_i\ge 0$ if and only if $1-a_i^\top \Omega a_i\ge 0$.  Property 3) follows directly from the fact that maximizing the log-determinant of $\Omega$ implies that the volume of the ellipsoid described by $\mathcal{B}$ is maximized. Finally, the last statement follows directly from properties 1) and 2), concluding the proof.
$\hfill \blacksquare$

\subsection{Proof of Lemma \ref{lemma:super:general}}
Our goal is to construct a non-negative supermartingale $\zeta_t$ starting from \eqref{eq:safety:martingale}. To this end, note that \eqref{eq:safety:martingale} is equivalent to
$$
\mathbb{E}[x_{t+1}^\top \Omega^{-1}x_{t+1}] \leq \tilde{\beta} (x_t^\top \Omega^{-1} x_t) + \psi,
$$
for $\tilde{\beta} = 1 - \beta \in (0,1)$ and $\psi =  \beta - \delta \geq 0$ since $\delta \in (\beta-1, \beta]$. By multiplying both sides of the previous inequality by $\prod_{i=1}^{t+1} g_i(\beta,\delta)$, which is strictly greater than zero by property i), and adding the term $\left( H - \sum_{i=1}^{t+1}h_i(\beta,\delta)\right)$, which is greater or equal than zero by definition of $H$, we get
$$
\begin{aligned}
& \mathbb{E}\left[ \left(\prod_{i=1}^{t+1}g_i(\beta,\delta) \right)(x_{t+1}^\top \Omega^{-1}x_{t+1}) + \left( H - \sum_{i=1}^{t+1}h_i(\beta,\delta)\right)\right] \\
\leq & \left(\prod_{i=1}^{t+1}g_i(\beta,\delta) \right) \left(\tilde{\beta} (x_t^\top \Omega^{-1} z_t) + \psi\right)+ \left( H - \sum_{i=1}^{t+1}h_i(\beta,\delta)\right)\\
\aleq & \left(\prod_{i=1}^{t}g_i(\beta,\delta) \right) (x_t^\top \Omega^{-1} x_t) +  \\
& \quad \quad \underbrace{\left(\prod_{i=1}^{t}g_i(\beta,\delta) \right) \psi + \left( H - \sum_{i=1}^{t+1}h_i(\beta,\delta)\right)}_{(\star)}\\
\end{aligned}
$$
where $(a)$ follows since $ g_{t+1}(\beta,\delta) \tilde{\beta} \leq 1$ by property i). Next, notice that 
$$
\begin{aligned}
&\left(H - \sum_{i=1}^t h_i(\beta,\delta)\right) - \left(H - \sum_{i=1}^{t+1}h_i(\beta,\delta) \right)\\
=&\sum_{i=1}^{t+1} h_i(\beta,\delta) - \sum_{i=1}^t h_i(\beta,\delta) =  h_{t+1}(\beta,\delta) \geq 0.
\end{aligned}
$$
Since $\prod_{i=1}^{t+1}g_i(\beta,\delta) \psi \leq h_{t+1}(\beta,\delta)$ by property ii), we conclude that 
$$
(\star) \leq \left( H - \sum_{i=1}^t h_i(\beta,\delta) \right).
$$
As a consequence, by definition we have $\mathbb{E}[\zeta_{t+1}|\zeta_t] \leq \zeta_t$, hence $\zeta_t$ is a supermartingale. The non-negative follows easily by properties i) and ii), the definition of $H$ and by observing that  $x_t^\top \Omega^{-1}x_t \geq 0$ for all $x_t \in \mathbb{R}^n$ since $\Omega^{-1} \succ 0$. Finally, note that $w_t$ is sub-Gaussian by Assumption \ref{ass:noise}; thus, $x_t$ is as well due to the closure of sub-Gaussian distributions under linear transformations. Then, we compactly rewrite $\zeta_t = G_t x_t^\top \Omega^{-1}x_t + \tilde{H}_t$, for $G_t > 0,\: \tilde{H}_t \geq 0$ by the properties i), ii) and the definition of $H$. We notice that $G_t, \tilde{H}_t$ do not change the tail behavior, thus we focus on establishing the finiteness of $\mathbb{E}[|x_t^\top \Omega^{-1}x_t|] = \mathbb{W}[x_t^\top \Omega^{-1}x_t]$ since $\Omega \succ 0$. The claim immediately follows the second moment of $x_t$ is finite under Assumption \ref{ass:noise}, hence $\mathbb{E}[|\zeta_t|] < \infty$.
$\hfill \blacksquare$

\subsection{Proof of Proposition \ref{prop:Ville:general}}
We shall exploit $\zeta_t$ constructed in Lemma \ref{lemma:super:general}. By Lemma \ref{lemma:Ville}), for any given $a > 0$ it holds
\begin{equation}\label{eq:fin:1}
\text{Pr}\left(\sup_{t \in \mathcal{T}} \zeta_t > a \right) \leq \mathbb{E}[\zeta_0]/a.
\end{equation}
Recall the safety description in \eqref{eq:safety:finite:recap}. Our goal is to then find an $a^\star$ such that 
\begin{equation}\label{eq:fin:1a}
\left(\sup_{t in \mathcal{T}} \zeta_t \leq  a^\star \right) \implies \left( \inf_{t \in \mathcal{T}} b(x_t) \geq -\gamma \right),
\end{equation}
where $\gamma > 0$ describes the super-level set representing the safe set. In our case, we fix $\gamma = 0$ (and hence drop it consequently); however, one can use the same argument to analyze any super-level set of interest.  Note that by \eqref{eq:fin:1} and \eqref{eq:fin:1a}, we can bound the safety probability as
\begin{equation}\label{eq:fin:2}
\begin{aligned}
 & \text{Pr}(x_t \in \mathcal{B}, \forall t \in \mathcal{T}) \\
  &  \quad = \text{Pr}\left( \inf_{t \in \mathcal{T}} b(x_t) \geq 0 \right) 
\geq   \text{Pr}\left(\sup_{t \in \mathcal{T}} \zeta_t \leq a^\star \right) \\
&\quad  =  1 - \text{Pr}\left(\sup_{t \in \mathcal{T}} \zeta_t > a^\star \right)
\geq  1 - \mathbb{E}[\zeta_0]/a^\star.
\end{aligned}
\end{equation}
\vspace{-0.3cm}
Consider 
\begin{equation}
a(t) = \left( \prod_{i=0}^t g_i(\delta, \beta) \right) + \left( H - \sum_{i=0}^t h_i(\beta, \delta) \right)
\end{equation}
and define $a^\star = \min_{t \in \mathcal{T}} a_t$. Since in \eqref{eq:fin:1a} we assume $\sup_{t \in \mathcal{T}} \zeta_t  \leq a^\star$ holds, we have
$$
\begin{aligned}
0 \geq & \zeta_t - a^\star \ageq  \zeta_t - a(t) \\
= & \left( \prod_{i=1}^{t+1}g_i(\beta,\delta) (x_t^\top \Omega^{-1} x_t) + H - \sum_{i=1}^t h_i(\beta,\delta)\right) \\
& -\left( \prod_{i=0}^t g_i(\delta, \beta) +  H - \sum_{i=0}^t h_i(\beta, \delta) \right)\\
= &\left( \prod_{i=1}^{t+1}g_i(\beta,\delta)\right) (x_t^\top \Omega^{-1} x_t - 1)\\
= & \left( \prod_{i=1}^{t+1}g_i(\beta,\delta)\right) \left(-b(x_t)\right),
\end{aligned}
$$
where (a) follows by definition of minimum. Since $g_i(\beta,\delta) > 0$ for all $i$ by property i) in Lemma \ref{lemma:super:general}, it holds $b(x_t) \geq 0, \: \forall t \in \mathcal{T}$ impliying $\inf_{t \in \mathcal{T}} b(x_t) \geq 0$. At this point, invoking Ville's Lemma on $\zeta_t$ constructed in Lemma \ref{lemma:super:general}, we get
$
 \text{Pr}(\sup_{t \in \mathcal{T}} \zeta_t > a^\star) \leq \alpha
$
and as a consequence
$
 \text{Pr}(x_t \in \mathcal{B}, \: \forall t \in \mathcal{T}) \geq 1 - \alpha,
$
where $\alpha$ is as in the proposition's statement. $\hfill \blacksquare$

\subsection{Proof of Theorem \ref{th:convex-design:pp}}
We first prove that \eqref{eq:convex-fin-n}, \eqref{eq:convex-fin-m} and \eqref{eq:convex-fin-e} provide a convex reformulation for \eqref{eq:safety:martingale}. By Assumption \eqref{ass:noise} and the Markovian property of our system, condition \eqref{eq:safety:martingale} under the parametrization in \eqref{eq:parametrization:all} evaluates as
$$
\begin{aligned}
&1 - (Ax + BY\Omega^{-1}x)^\top \Omega^{-1}(Ax + BY\Omega^{-1}x) \\
& - \tilde{\beta}(1 - x^\top \Omega^{-1}x)- \delta - \text{Tr}(\Omega^{-1} \Sigma) \geq 0
\end{aligned}
$$
where $\tilde{\beta} = 1 - \beta \in (0,1)$. This is equivalent to
$$
\begin{cases}
 x^\top (\tilde{\beta} \Omega^{-1} - (A + BY\Omega^{-1})^\top \Omega^{-1}(A+BY\Omega^{-1})) x \succeq 0\\
 \beta - \delta - \text{Tr}(\Omega^{-1}\Sigma) \geq 0.
\end{cases}
$$
We focus on the first condition. From \eqref{eq:convex-fin-b}, we directly have $\Omega$ is invertible. Pre- and post-multiplying the matrix by $\mathrm{diag}(\Omega,I)$, we get
    \begin{align*}
        &\begin{bmatrix}
            \tilde{\beta}\Omega^{-1}&A^\top \Omega^{-1}+\Omega^{-1}B^\top Y^\top \Omega^{-1}\\
            \star&\Omega^{-1}
        \end{bmatrix}\succeq 0,\forall x\in\mathbb{R}^n\\
        &\begin{bmatrix}
            \tilde{\beta}\Omega^{-1}&(A+BY\Omega^{-1})^\top \Omega^{-1}\\
            \star&\Omega^{-1}
        \end{bmatrix}\succeq 0,\forall x\in\mathbb{R}^n.
    \end{align*}
    Applying a congruent transformation $\mathrm{diag}(I_n,\Omega)$, we get
    \begin{equation*}
        \begin{bmatrix}
            \tilde{\beta}\Omega^{-1}&(A+BY\Omega^{-1})^\top \Omega^{-1}\\
            \star&\Omega
        \end{bmatrix}\succeq 0,\forall x\in\mathbb{R}^n.
    \end{equation*}
    Given that $\beta>0$ and $\Omega\succ 0$, we can apply Schur complement to obtain an equivalent condition:
    \begin{equation}
        \tilde{\beta}\Omega^{-1}-\star \Omega^{-1}(A+BY\Omega^{-1})\succeq 0,\forall x\in\mathbb{R}^n.
    \end{equation}
    By multiplying vector $x$ on both sides of the matrix, we conclude the first case. As for the second case, by Schur complement it equals to
    $$
    \beta - \delta - \lambda \geq 0, \: \begin{bmatrix}\lambda I & \Sigma^{1/2} \\ \star & \Omega \end{bmatrix} \succeq 0
    $$
Next, we encode the condition $\mathcal{I} \subseteq \mathcal{B}$. Note that we can equivalently write $\mathcal{I} = \{1 - x^\top \frac{R}{1-\sigma} x \geq 0\}$ with $\frac{R}{1-\sigma} \succ 0$ since $\sigma \in (0,1)$. The set inclusion then reads
\begin{equation}\label{eq:initial}
\{ x \in \mathbb{R}^n | 1 - x^\top \frac{R}{1-\sigma}  x\} \subseteq \{ x \in \mathbb{R}^n | 1 - x^\top \Omega^{-1} x\},
\end{equation}
which is trivially satisfied as \eqref{eq:convex-fin-f} holds.
Finally, $\mathcal{B} \subseteq \mathcal{S}$ is encoded via \eqref{eq:convex-fin-h} by exploiting Farkas' lemma, as done in Theorem \eqref{th:convex-design:pp}. The proof concludes by invoking Proposition \eqref{prop:Ville:general} together with the design rules in Example 1 and noticing that $\mathcal{B} \subseteq \mathcal{S}$ implies $\text{Pr}(x_t \in \mathcal{B}, \: \forall t \in \mathcal{T}) \geq 1 - \alpha.$
 $\hfill \blacksquare$

\section*{Proofs of Section 4}
\subsection{Proof of Lemma \ref{lemma:input:constraint}}
 Eq.\eqref{eq:convex-lin-input} is sufficient to
    $        \mathcal{H}_{11}^i-\mathcal{H}_{12}I_{n+1}\mathcal{H}_{12}^\top \succeq 0,i=1,\ldots,k,\forall x\in\mathbb{R}^n,
   $
    which implies
    \begin{equation*}
        \Omega-\frac{Y^\top H_i^\top H_iY }{2\lambda h_i-\lambda^2}\succeq 0,i=1,\ldots,k,\forall x\in\mathbb{R}^n.
    \end{equation*}
    Given that $\lambda >0$, we have
    \begin{equation*}
        \lambda\Omega-\frac{Y^\top H_i^\top H_iY}{2h_i-\lambda}\succeq 0,i=1,\ldots,k,\forall x\in\mathbb{R}^n.
    \end{equation*}
    Left- and right- multiply the matrix by $\Omega^{-1}$, and define $K=Y\Omega^{-1}$, we obtain
    \begin{equation*}
        \lambda\Omega^{-1}-\frac{K^\top H_i^\top H_iK}{2h_i-\lambda}\succeq 0,i=1,\ldots,k,\forall x\in\mathbb{R}^n.
    \end{equation*}
    Apply Schur complement again, we have for all $x\in\mathbb{R}^n$
    \begin{align*}
        &\begin{bmatrix}
            \lambda\Omega^{-1}&-K^\top H_i^\top \\-H_iK&2h_i-\lambda
        \end{bmatrix}\succeq 0,i=1,\ldots,k,\forall x\in\mathbb{R}^n \Longleftrightarrow\\
       &-2(H_iKx-h_i)+\lambda(x^\top \Omega^{-1}x-1)\ge 0,i=1,\ldots,k..
    \end{align*}
    Then we have for $i=1,\ldots,k$, $H_iu(x)\le h_i$ for any $x$ such that $b(x)\ge 0$. Hence, we conclude the proof. 
$\hfill \blacksquare$

\subsection{Proof of Proposition \ref{prop:gelbrich}}
We first rewrite \eqref{eq:robust:cond} in a more convenient, yet equivalent, manner. To this end, fix a $x \in \mathbb{R}^n$ and let $M = x^\top (\tilde{\beta} \Omega^{-1} - \star \Omega^{-1}(A+BY\Omega^{-1}))x + (1-\tilde{\beta} + \delta)$. Next, define $$
f(w) \coloneqq \begin{bmatrix} w \\ 1\end{bmatrix}^\top \begin{bmatrix} \Omega^{-1} & 0 \\ 0 & -M \end{bmatrix} \begin{bmatrix} w \\ 1\end{bmatrix}. $$
It is easy to see that \eqref{eq:robust:cond} is equivalent to 
\begin{equation}\label{eq:equiv:cond}
\sup_{\nu \in \mathcal{A}_\rho(\hat{\mu})} \mathbb{E}_{w \sim \nu} f(w) \leq 0.
\end{equation}
By \cite[Theoreom~6]{nguyen2021mean}, \eqref{eq:equiv:cond} admits the strong dual
$$
\begin{cases}
\inf & q_0+\gamma\left(\rho^2-\operatorname{Tr}[S]\right)+z+\operatorname{Tr}[Z] \\
\text { s.t. } & \gamma \in \mathbb{R}_{+}, z \in \mathbb{R}_+, Z \in \mathbb{S}_{+}^d,\left(q_0, Q\right) \in \mathcal{Q} \\
& \left[\begin{array}{cc}
 \gamma I-Q & \gamma S^{\frac{1}{2}} \\
 \gamma S^{\frac{1}{2}} & Z
 \end{array}\right] \succeq 0, 
\end{cases}
$$
where $\mathcal{Q} = \{ q_0 \in \mathbb{R}, Q \in \mathbb{S}^d : q_0 +  w^\top Q w \geq f(w), \forall w \in \mathbb{R}^d\}$. Exploiting the structure of $f(w)$, the constraint $(q_0, Q) \in \mathcal{Q}$ becomes
$$
\begin{bmatrix} w \\ 1\end{bmatrix}^\top \begin{bmatrix} Q & 0 \\ 0 & q_0 \end{bmatrix} \begin{bmatrix} w \\ 1\end{bmatrix} \geq \begin{bmatrix} w \\ 1\end{bmatrix}^\top \begin{bmatrix} \Omega^{-1} & 0 \\ 0 & -M \end{bmatrix} \begin{bmatrix} w \\ 1\end{bmatrix},
$$
that is equivalent to $Q \succeq \Omega^{-1}$ and $q_0 \succeq -M$. The first inequality can be reformulated by Schur complement as $$\begin{bmatrix} Q & I \\ I & \Omega \end{bmatrix} \succeq 0.$$
As for the second inequality, we shall enforce it for all $x \in \mathbb{R}^n$. To this end, following a similar reasoning as in the proof of Theorem 3.7 yields $$q_0I \succeq \begin{bmatrix}
        \tilde{\beta}\Omega & \Omega A^\top+B^\top Y^\top \\
        \star&\Omega
    \end{bmatrix},$$ concluding the proof.  $\hfill \blacksquare$
    
We additionally reason on the connection between the robustification provided by the ambiguity set \eqref{eq:ambiguity:set} based on the Gelbrich distance and other types of robustification. Specifically, by \cite[Proposition~3]{masarotto2019procrustes} the following relationships holds true for any $A_1, A_2 \in \mathbb{S}^d_+$
\begin{equation}\label{eq:relat}
\begin{aligned}
\sqrt{\text{Tr}\left[A_1 + A_2 - 2(A_1^{1/2}A_2A_2^{1/2})^{1/2}\right]}\leq \|A_1^{1/2} -  A_2^{1/2}\|_F, 
\end{aligned}
\end{equation}
where $\|\cdot\|_F$ denotes the Frobenius norm.
As a result, one can alternatively require the forward invariant condition in \eqref{eq:safety:martingale} to hold for all $\Sigma \succ 0\; : \; \|\Sigma^{1/2} - S^{1/2}\|_F \leq \rho$. Based on \eqref{eq:relat}, the latter implies \eqref{eq:robust:cond} with the Gelbrich ambiguity set in \eqref{eq:ambiguity:set} of radius $\rho$. 

Note that this robust constraint requires the following implication to hold 
\begin{equation}\label{eq:robust:1}
\|\Sigma^{1/2} - S^{1/2}\|_F \leq \rho \implies  \eqref{eq:safety:martingale}.
\end{equation}
This implication can be enforced following a similar procedure as in the proof of Theorem 3.7 and invoking the S-Procedure (Lemma \ref{eq:S-procedure1}) to show that 
\begin{equation}\label{eq:robust:2}
\|\Sigma^{1/2} - S^{1/2}\|_F \leq \rho \implies  \begin{bmatrix}
\lambda I & \Sigma^{ 1/2} \\ \star & \Omega
\end{bmatrix} \geq 0
\end{equation}
 holds if and only if there exists a $\nu > 0$ such that
\begin{equation}\label{eq:robust:sigma}
 \begin{bmatrix}
(\lambda - \nu) I & S \\ S^\top & \Omega - \nu I
\end{bmatrix} \geq 0, \quad   \nu \geq \frac{\rho^2}{2}.
\end{equation}
Note that \eqref{eq:robust:sigma} consists of a set of convex constraints in $(\lambda, \nu, \Omega)$ that can be directly embedded into Problem \eqref{eq:convex-without-input-constraints-finite}.

\end{document}